\documentclass[leqno,12pt]{article}
\usepackage{amsmath,amssymb,rotating,a4wide,color}
\usepackage[all,cmtip]{xy}

\usepackage{verbatim}

\NoComputerModernTips

\newdir^{ (}{{}*!/-3pt/\dir^{(}}
\newdir_{ (}{{}*!/-3pt/\dir_{(}}

\entrymodifiers={+!!<0pt,\fontdimen22\textfont2>}

\newbox\circbulletbox
\setbox\circbulletbox=\hbox{{\large$\circledcirc$}\kern-8.7pt\raise .6pt\hbox{$\bullet$}\,}
\def\circbullet{\copy\circbulletbox}

\def\mysmash#1#2{\smash{\raisebox{#1}{$\scriptstyle #2$}}}

\def\mybullet{\makebox[0pt]{\raisebox{0pt}[0pt][0pt]{$\bullet$}}}

\let\le\leqslant
\let\ge\geqslant

\def\mycirc{{\kern1pt\circ\kern2pt}}

\def\half{{\frac{1}{2}}}
\def\sm{{\rm sm}}
\def\barX{{\overline{X}}}
\def\gensigma{\langle\sigma\rangle}
\def\image{\mathop{\rm image}\nolimits}

\def\Aut{\mathop{\rm Aut}\nolimits}

\def\Spec{\mathop{\rm Spec}\nolimits}

\def\mod{\mathop{\rm mod}\nolimits}

\def\diag{\mathop{\rm diag}\nolimits}
\def\ord{\mathop{\rm ord}\nolimits}

\def\PGL{\mathop{\rm PGL}\nolimits}

\def\et{{\rm \acute{e}t}}

\def\id{{\rm id}}

\let\phi\varphi
\let\theta\vartheta
\let\epsilon\varepsilon
\let\setminus\smallsetminus

\newtheorem{Thm}{Theorem}[section]
\newtheorem{Prop}[Thm]{Proposition}
\newtheorem{Lem}[Thm]{Lemma}

\newtheorem{Cor}[Thm]{Corollary}

\newtheorem{Def}[Thm]{Definition}

\newtheorem{Rem}[Thm]{Remark}
\newtheorem{Ex}[Thm]{Example}

\newtheorem{Ques}[Thm]{Question}
\newtheorem{Class}[Thm]{Classification}

\def\UseTheoremCounterForNextEquation{\setcounter{equation}{\value{Thm}}\addtocounter{Thm}{1}}

\def\qed{{\hskip0pt\unskip\unskip\nobreak\hfil\penalty50
          \hskip1em\hbox{}\nobreak\hfil
           {$\square$}
          \parfillskip=0pt\finalhyphendemerits=0
          \par}\medskip}

\newenvironment{Proof}
               {\noindent{\bf Proof.}\ }
               {\qed}

               {\noindent{\bf Proof of #1.}\ }
               {\qed}

\newcommand{\BA}{{\mathbb{A}}}

\newcommand{\BC}{{\mathbb{C}}}

\newcommand{\BF}{{\mathbb{F}}}
\newcommand{\BG}{{\mathbb{G}}}

\newcommand{\BP}{{\mathbb{P}}}
\newcommand{\BQ}{{\mathbb{Q}}}

\newcommand{\BZ}{{\mathbb{Z}}}

\newcommand{\CO}{{\cal O}}

\newbox\mybox
\def\arrover#1{\mathrel{
       \setbox\mybox=\hbox spread 1.4em
              {\hfil$\scriptstyle#1$\hfil}
       \vbox{\offinterlineskip\copy\mybox
             \hbox to\wd\mybox{\rightarrowfill}}}}

\def\larrover#1{\mathrel{
       \setbox\mybox=\hbox spread 1.4em
              {\hfil$\scriptstyle#1\vphantom{g}$\hfil}
       \vbox{\offinterlineskip\copy\mybox
             \hbox to\wd\mybox{\leftarrowfill}}}}

\def\ontoover#1{\mathrel{
       \setbox\mybox=\hbox spread 1.4em
              {\hfil$\scriptstyle#1\vphantom{g}$\hfil}
       \vbox{\offinterlineskip\copy\mybox
             \hbox to\wd\mybox{\rightarrowfill\hskip-2.8mm
                               $\rightarrow$}}}}
\def\leftontoover#1{\mathrel{
       \setbox\mybox=\hbox spread 1.4em
              {\hfil$\scriptstyle#1\vphantom{g}$\hfil}
       \vbox{\offinterlineskip\copy\mybox
             \hbox to\wd\mybox{$\leftarrow$\hskip-2.8mm
                               \leftarrowfill}}}}
\let\longto\longrightarrow
\let\into\hookrightarrow
\let\onto\twoheadrightarrow

\def\Bigskip{\bigskip\bigskip}

\begin{document}

\title{Finiteness and liftability of postcritically finite 
quadratic morphisms in arbitrary characteristic}
%\title{All postcritically finite quadratic morphisms \\
%can be lifted to characteristic zero}

\author{Richard Pink\\[12pt]
\small Department of Mathematics \\[-3pt]
\small ETH Z\"urich\\[-3pt]
\small 8092 Z\"urich\\[-3pt]
\small Switzerland \\[-3pt]
\small pink@math.ethz.ch\\[12pt]}

% \footnotetext[1]{Dept. of Mathematics, ETH Z\"urich, 8092 Z\"urich, Switzerland, {\tt pink@math.ethz.ch}}

\date{August 26, 2013}
\maketitle

\begin{abstract}
We show that for any integer $n$ and any field $k$ of characteristic $\not=2$ there are at most finitely many isomorphism classes of quadratic morphisms from $\BP^1_k$ to itself with a finite postcritical orbit of size~$n$. As a consequence we prove that every postcritically finite quadratic morphism over a field of positive characteristic can be lifted to characteristic zero with the same combinatorial type of postcritical orbit. The associated profinite geometric monodromy group is therefore the same as in characteristic zero, where it can be described explicitly by generators as a self-similar group acting on a regular rooted binary tree.
\end{abstract}
%
%\newpage
%{\advance\baselineskip by -6pt
%\tableofcontents
%}

{\renewcommand{\thefootnote}{}
}

%%%%%%%%%%%%%%%%%%%%%%%%%%%%%%%%%%%%%%%%%%%%%%%%%%%%%%%%%%%%%%%%%%%%%%%%%%%%%%%%%%%%%%%%%%%%%

\newpage
\addtocounter{section}{-1}
\section{Introduction}
\label{Intro}

The dynamics of postcritically finite quadratic polynomials and quadratic rational maps on $\BP^1$ have been studied from various angles. This author is interested in arithmetic properties of the associated monodromy representations. In modern arithmetic geometry one tries to understand a situation in characteristic zero and in positive characteristic on the same footing and expects that each of them sheds light on the other. The present paper is meant to be a contribution in this direction. 

\medskip
As a consequence of Thurston rigidity (see Douady-Hubbard \cite{DouadyHubbard1993}, Brezin et al. \cite{BBLPP2000}) it is known that for any integer $n$ there are at most finitely many isomorphism classes of quadratic rational maps $\BP^1_\BC\to\BP^1_\BC$ with a postcritical orbit of size~$n$. Our main result Theorem \ref{GammaModuliQuasiFinite} can be phrased as saying that the same is true over any algebraically closed field $k$ of characteristic $\not=2$.  After finishing the work on this paper, we learned that this was already proved using different methods by Benedetto-Ingram-Jones-Levy \cite[Cor.~6.3]{BIJL2012}, see also Levy \cite{Levy2012}.

We prove Theorem \ref{GammaModuliQuasiFinite} in a purely algebro-geometric fashion, though the proof is surprisingly complicated. Assume that the finiteness is not true. Then an easy argument involving moduli spaces shows that there must exist a non-isotrivial family of quadratic morphisms (we prefer the precise terminology of algebraic geometry over the unspecific word `map') with postcritical orbits of the same combinatorial type over some smooth curve $D$ over~$k$. As the moduli space is affine, this family must have bad reduction at some point of a smooth compactification $\bar D$ of~$D$. In fact we exhibit a point where the reduction has a certain combinatorial type (see the proof of Lemma \ref{CrossRatio}), and by a different argument we show that such a type is actually impossible, thus arriving at a contradiction.

The combinatorial properties of the reduction are analyzed by means of stable marked curves. For this observe that the postcritical orbit defines a finite collection of disjoint sections of $\BP_D^1$ over~$D$, turning $\BP^1_D$ into a smooth marked curve of genus zero over~$D$. As such, it possesses a unique extension to a stable marked curve with disjoint sections over~$\bar D$. The degenerate fibers are trees of rational curves among which the marked points are distributed in a certain fashion, whose combinatorics can be described by what may be called stable marked trees. Although the given quadratic morphism $\BP^1_D\to\BP^1_D$ does not extend to the stable curve, it extends as a correspondence using a certain augmented stable marked extension, with the help of which we describe the precise combinatorial effect of the quadratic morphism on the respective marked trees (see Section \ref{StableCurvesQuad}).
% My $\Gamma$-marked trees with the maps $T\into\tilde T\onto T$ 
(There may be a relation with the mapping trees from Pilgrim \cite[\S2.1]{Pilgrim2003}, but we have not made a detailed comparison.)
% There may also be a vague relation with Hubbard trees (see Poirier 1993), although Hubbard trees have angles and their set of vertices is equal to~$\Gamma$.
This description suffices to exclude certain types of bad reduction by purely combinatorial arguments (see Section \ref{Exclude}). At some step we also use the local multiplicity of a singular point on the total space of the family, which is not only a combinatorial invariant of the special fiber (see Proposition \ref{QuotientProp} (b) and the proof of Lemma \ref{B'EqualA}).
All these arguments make up the bulk of this paper from Section \ref{Finite} onwards and serve only to establish Theorem  \ref{GammaModuliQuasiFinite}.

The moduli space that we use is a rigidified version of the moduli space of dynamical systems of degree $2$ from Silverman \cite[\S4.6]{Silverman1986}. Using it we construct a moduli space $M_\Gamma$ of $\Gamma$-marked quadratic morphisms for any finite mapping scheme $\Gamma$ of cardinality~${\ge3}$. The points of $M_\Gamma$ over an algebraically closed field $k$ of characteristic $\not=2$ are in bijection with the isomorphism classes of quadratic morphisms over $k$ whose postcritical orbit is combinatorially equivalent to~$\Gamma$. Our finiteness theorem \ref{GammaModuliQuasiFinite} then becomes the statement that all fibers of $M_\Gamma$ over $\Spec\BZ[\half]$ are finite.

On the other hand, it is not hard to construct $M_\Gamma$ as the joint zero locus of two polynomials in a suitable Zariski open subscheme of $\BA^2\times\Spec\BZ[\half]$ (see Proposition \ref{GammaModuli}).  Using the local flatness criterion from commutative algebra this together with the finiteness implies that $M_\Gamma$ is flat over $\Spec\BZ[\half]$. As a consequence, any $\Gamma$-marked quadratic morphism over a field of positive characteristic can be lifted to characteristic zero (see Corollary \ref{GammaModuliLift}).

Finally, consider a quadratic morphism $\BP^1_k\to\BP^1_k$ over an algebraically closed field $k$ of characteristic $>2$ with finite postcritical orbit $\Sigma\subset\BP^1_k$. The associated monodromy representation is a homomorphism from the (geometric) \'etale fundamental group ${\pi_{1,\et}(\BP^1_k\setminus \Sigma)}$ to the automorphism group of a regular rooted binary tree~$T$ (see Section \ref{MonodromyGroups}). By combining the liftability with Grothendieck's theorem on the specialization of the tame fundamental group we show that the image of this homomorphism is the same as for some quadratic morphism with a combinatorially equivalent postcritical orbit over~$\BC$. The latter is simply the closure in the profinite topology of $\Aut(T)$ of the image of the usual fundamental group $\pi_1(\BP^1(\BC)\setminus \Sigma)$ and can be described by explicit generators as a self-similar subgroup of $\Aut(T)$ as, say, in 
Bartholdi-Nekrashevych \cite{Bartholdi-Nekrashevych-2008},
Dau \cite{Dau}, 
Grigorchuk et al.~\cite{Grigorchuk-et-al-2007},
Nekrashevych \cite[Ch.5]{Nekrashevych-2005}, \cite{Nekrashevych-2009}.

\medskip
There are many interesting open questions that one might pursue next. First, while almost all combinatorially conceivable finite mapping schemes seem to occur as the postcritical orbit of some quadratic morphism over~$\BC$, some of them cannot occur in all characteristics (see Remarks \ref{AllOccur} and \ref{NotAllOccur}). So which mapping schemes occur in which characteristic, and why?

Second, some $\Gamma$-marked quadratic morphisms over a field possess non-trivial infinitesimal deformations (see Remark \ref{NotReduced}). When precisely does this occur? What are the conceptual reasons for it?

Third, we have analyzed the possible combinatorial types of degeneration of $\Gamma$-marked quadratic morphisms only to the extent necessary to prove the finiteness theorem. We have excluded certain types, while others are possible, as in Remark \ref{NotAllOccur} and Example \ref{NotAllOccur2}. So which kinds of degenerations actually occur? 

Fourth, how is the combinatorial type of degeneration over a non-archimedean local field, as described in the present paper, related to the dynamical properties of the quadratic morphism in the generic fiber?

Fifth, the above mentioned consequence for the monodromy group concerns the geometric fundamental group, i.e., the case where the base field $k$ is algebraically closed. What can be said about the image of the arithmetic fundamental group, i.e., when $k$ is not algebraically closed? What about the images of Frobenius elements?

Sixth, how can the method be generalized to morphisms of higher degree?

%%%%%%%%%%%%%%%%%%%%%%%%%%%%%%%%%%%%%%%%%%%%%%%%%%%%%%%%%%%%%%%%%%%%%%%%%%%%%%%%%%%%%%%%%%%%%

\Bigskip
\noindent{\bf\Large Part I: Moduli of quadratic morphisms}

%%%%%%%%%%%%%%%%%%%%%%%%%%%%%%%%%%%%%%%%%%%%%%%%%%%%%%%%%%%%%%%%%%%%%%%%%%%%%%%%%%%%%%%%%%%%%

\section{Stable quadratic morphisms}
\label{StableQM}

Let $S$ be a scheme over $\Spec\BZ[\half]$.

\begin{Def}\label{QuadDef}
A \emph{quadratic morphism (with marked critical points)} over $S$ is a quadruple $(C,f,P,Q)$ where
\begin{enumerate}
\item[(a)] $C$ is a curve over $S$ that is locally for the Zariski topology on $S$ isomorphic to $\BP^1\times S$,
\item[(b)] $f\colon C\to C$ is a morphism over $S$ which is fiberwise of degree~$2$, and
\item[(c)] $P$, $Q\in C(S)$ are sections whose images are precisely the critical points of~$f$, that is, the points where $df=0$.
\end{enumerate}
\end{Def}

\begin{Rem}\label{QuadDefRem}
\rm In Milnor \cite[\S6]{Milnor1993} this is called a `critically marked quadratic rational map'. In the present paper, all quadratic morphisms are endowed with marked critical points, but for brevity we will simply speak of quadratic morphisms. We will often speak of a quadratic morphism $f$ if the other data $C$, $P$, $Q$ are not explicitly used or are understood.
\end{Rem}

There is an obvious notion of isomorphisms and automorphisms of quadratic morphisms with marked critical points.
Note that the sections $P$, $Q$ in any quadratic morphism are fiberwise distinct, because $f$ possesses precisely two critical points in every geometric fiber. 
For any integer $n\ge0$ we let $f^n$ denote the $n^{\rm th}$ \emph{iterate of~$f$}, obtained by setting $f^0:=\id$ and $f^{n+1}:=f\circ f^n$.
For any section $R\in C(S)$ we abbreviate $f^n(R) := f^n\circ R$. 
If the quadratic morphism is defined over field~$K$, the sections $P$, $Q$ are really $K$-valued points of~$C$, and so are their images $f^n(P)$, $f^n(Q)$ under all iterates of~$f$.

\begin{Def}\label{PostCritDef}
For any quadratic morphism $f$ over a field, the set ${\{ f^n(P), f^n(Q) \mid n\ge1 \}}$ is called the \emph{(strictly) postcritical orbit} of~$f$. If this set is finite, then $f$ is called \emph{post\-criti\-cally finite}.
\end{Def}

\begin{Prop}\label{Unstable}
For any quadratic morphism $(C,f,P,Q)$ over a field $K$ the following are equivalent:
\begin{enumerate}
\item[(a)] The postcritical orbit has cardinality $\ge3$.
\item[(b)] At least one of $f(P)$, $f(Q)$ is distinct from both $P$, $Q$.
\item[(c)] $(C,f,P,Q)$ is not isomorphic to $(\BP^1_K, x\mapsto ax^{\pm2},0,\infty)$ for any sign and any $a\in K^\times$.
\end{enumerate}
Moreover, these conditions imply:
\begin{enumerate}
\item[(d)] $\Aut(C,f,P,Q)=1$.
\end{enumerate}
\end{Prop}

\begin{Proof}
The implication (a)$\Rightarrow$(b) is trivial. To prove the converse we first note that since $f^{-1}(P)=\{P\}$ and $f^{-1}(Q)=\{Q\}$, for any point $R\in C(K)$ we have
\UseTheoremCounterForNextEquation
\begin{equation}\label{RamLift}
\biggl\{\begin{array}{l}
f(R)=f(P)\ \Longleftrightarrow\ R=P,\ \ \ \hbox{and} \\[3pt]
f(R)=f(Q)\ \Longleftrightarrow\ R=Q.
\end{array}\biggr\}
\end{equation}
In particular, since $P\not=Q$, we always have 
\UseTheoremCounterForNextEquation
\begin{equation}\label{Jamiri}
f(P)\not=f(Q).
\end{equation}
Now assume that (b) does not hold, i.e., that the (strictly) postcritical orbit has cardinality $\le2$. Then both $f^2(P)$ and $f^2(Q)$ must be among $f(P)$ and $f(Q)$. But by (\ref{RamLift}) this implies that both $f(P)$ and $f(Q)$ are among $P$ and~$Q$, and so (a) does not hold. Thus (a) and (b) are equivalent.
Furthermore, if (b) does not hold, the preceding remarks imply that $\{f(P),f(Q)\} = \{P,Q\}$. Identify $C$ with $\BP^1_K$ such that $P=0$ and $Q=\infty$. Then a direct calculation shows that $f$ must have the form excluded in (c). Conversely, the quadratic morphism excluded in (c) visibly satisfies $\{f(0),f(\infty)\} = \{0,\infty\}$. This proves the equivalence (b)$\Leftrightarrow$(c). 

Finally, any automorphism of $(C,f,P,Q)$ fixes each point in the postcritical orbit. Since any automorphism of $\BP^1_K$ fixing $\ge3$ points is trivial, we have  (a)$\Rightarrow$(d) and are done.
\end{Proof}

\begin{Def}\label{StableDef}
A quadratic morphism over $S$ is called \emph{stable} if in every fiber the postcritical orbit has cardinality $\ge3$.
\end{Def}

\begin{Prop}\label{StableModuli}
There is a fine moduli scheme $M$ of stable quadratic morphisms, and it is smooth of finite type of relative dimension $2$ over $\Spec\BZ[\half]$.
\end{Prop}

\begin{Proof}
First consider a quadratic morphism $(C,f,P,Q)$ over $S$ for which $f(Q)$ is disjoint from $P$ and~$Q$. Then there is a unique isomorphism $C\cong {\BP^1\times S}$ which sends the sections $(P,Q,f(Q))$ to $(0,\infty,1)$. After carrying out this identification, the quadratic morphism has the form $f(x) = \frac{x^2+a}{x^2+b}$ for unique $a$, $b\in\Gamma(S,\CO_S)$. Thus the isomorphism class of $(C,f,P,Q)$ is determined by $(a,b)$. Also, for any $a$, $b\in\Gamma(S,\CO_S)$ the formula $\frac{x^2+a}{x^2+b}$ defines a morphism $\BP^1\times S \to \BP^1\times S$ which is fiberwise of degree $2$ if and only if $a\not=b$ everywhere. Thus $(a,b)$ corresponds to a morphism $S\to\BA^2\setminus\diag(\BA^1)$, and conversely, to any such morphism we can associate the stable quadratic morphism $(\BP^1\times S, x\mapsto\frac{x^2+a}{x^2+b},0,\infty)$ with $f(\infty)=1$. This shows that the subfunctor of all stable quadratic morphisms satisfying $f(Q)\not=P,Q$ possesses a fine moduli scheme $M_2$ isomorphic to $(\BA^2\setminus\diag(\BA^1))\times\Spec\BZ[\half]$.

Next consider a quadratic morphism $(C,f,P,Q)$ over $S$ for which $f(P)$ is disjoint from $P$ and~$Q$. Then we identify $C$ with $\BP^1\times S$ by sending the sections $(P,Q,f(P))$ to $(0,\infty,1)$. Afterwards the quadratic morphism has the form $f(x) = \frac{cx^2+1}{dx^2+1}$ for unique $c$, $d\in\Gamma(S,\CO_S)$. Again we find that $(c,d)$ must represent a morphism $S\to\BA^2\setminus\diag(\BA^1)$, and conversely, that any such morphism determines a stable quadratic morphism $(\BP^1\times S, x\mapsto\frac{cx^2+1}{dx^2+1},0,\infty)$ with $f(0)=1$. This shows that the subfunctor of all stable quadratic morphisms satisfying $f(P)\not=P,Q$ possesses a fine moduli scheme $M_1$ which is also isomorphic to $(\BA^2\setminus\diag(\BA^1))\times\Spec\BZ[\half]$.

Now let $M_{21}\subset M_2$ and $M_{12}\subset M_1$ denote the open subschemes where the sections $f(P)$ and $f(Q)$ of the respective universal family of quadratic morphisms are both disjoint from $P$ and~$Q$. Then $M_{21}$ and $M_{12}$ represent the same functor and are therefore canonically isomorphic. Let $M$ be the scheme over $\Spec\BZ[\half]$ obtained by gluing $M_2$ and $M_1$ along this isomorphism. (An explicit calculation, not necessary for the proof, shows that the gluing isomorphism is given by 
$$\xymatrix@R-10pt@C-10pt{
M_1 \ar@{}[d]|{\wr\Vert\ } \ar@{^{ (}->}[r]
& M_{12} \ar@{}[d]|{\wr\Vert\ } \ar[rrr]^\sim
&&& M_{21} \ar@{}[d]|{\wr\Vert\ } & M_2 \ar@{}[d]|{\wr\Vert\ } \ar@{_{ (}->}[l] \\
\BA^2\setminus\diag(\BA^1) \ar@{}[r]|-{\supset} &
\BG_m^2\setminus\diag(\BG_m^1) 
\ar[rrr]^{(c,d) \mapsto (\frac{d^2}{c^3},\frac{d}{c^2})} 
&&& \BG_m^2\setminus\diag(\BG_m^1) \ar@{}[r]|-{\subset}  &
\BA^2\setminus\diag(\BA^1) \\}$$
where for brevity we have dropped the factor $\Spec\BZ[\half]$.)

Consider an arbitrary stable quadratic morphism $(C,f,P,Q)$ over $S$. Let $S_2$ be the open subscheme of~$S$ where $f(Q)$ is disjoint from $P$ and~$Q$, and $S_1$ the open subscheme of~$S$ where $f(P)$ is disjoint from $P$ and~$Q$. Then Proposition \ref{Unstable} implies that $S=S_1\cup S_2$. Moreover, the restriction of the family to $S_\nu$ is classified by a morphism $S_\nu\to M_\nu$, which induces a morphism from $S_1\cap S_2$ to $M_{12}\subset M_1$, respectively to $M_{21}\subset M_2$, and these morphisms are compatible with the given isomorphism $M_{12}\cong M_{21}$. Thus the morphisms $S_\nu\to M_\nu$ combine to a morphism $S\to M$, and $M$ is a fine moduli space for the moduli problem at hand.

Finally, the stated properties of $M$ result from the open covering by two copies of $(\BA^2\setminus\diag(\BA^1))\times\Spec\BZ[\half]$. 
\end{Proof}

%%%%%%%%%%%%%%%%%%%%%%%%%%%%%%%%%%%%%%%%%%%%%%%%%%%%%%%%%%%%%%%%%%%%%%%%%%%%%%%%%%%%%%%%%%%%%

\section{Mapping schemes}
\label{MappingSchemes}

\begin{Def}\label{AbsPostCritDef}
A \emph{finite mapping scheme} is a quadruple $(\Gamma,\tau,i_1,j_1)$ consisting of a finite set~$\Gamma$, a map $\tau\colon\Gamma\to\Gamma$ and two distinct elements $i_1$, $j_1\in\Gamma$ such that, abbreviating throughout $i_n := \tau^{n-1}(i_1)$ and $j_n := \tau^{n-1}(j_1)$ for all integers $n\ge2$, we have:
\begin{enumerate}
\item[(a)] $\Gamma = \{ i_n,\;j_n \mid n\ge1 \}$.
\item[(b)] For all $\gamma\in\Gamma$ we have $|\tau^{-1}(\gamma)| \le 2$.
\item[(c)] $|\tau^{-1}(i_1)|\le1$ and $|\tau^{-1}(j_1)|\le1$.
\end{enumerate}
\end{Def}

This is essentially the special case of degree $2$ of the definition in Brezin et al. \cite[Def.\,2.3]{BBLPP2000}, but the definitions in the literature vary. Often we will speak of a mapping scheme $\Gamma$ if the other data $\tau$, $i_1$, $j_1$ is understood.

\begin{Prop}\label{PostCritIsAbs}
The postcritical orbit of any postcritically quadratic morphism $f$ over a field, with the map induced by $f$ and the distinguished elements $i_1:=f(P)$ and $j_1:=f(Q)$, is a finite mapping scheme in the sense of Definition \ref{AbsPostCritDef}.
\end{Prop}

\begin{Proof}
By (\ref{Jamiri}) we have $f(P)\not=f(Q)$, and (a) holds by Definition \ref{PostCritDef}. Condition (b) follows from the fact that $f$ has degree~$2$, and (c) follows from (\ref{RamLift}).
\end{Proof}

\begin{Class}\label{AbsPostCritClass}
\rm The isomorphism classes of finite mapping schemes fall into different types according to which of the elements $i_n$ and $j_n$ are equal. Note that any relation of the form $i_n=i_m$ implies $i_{n+\ell}=i_{m+\ell}$ for all $\ell\ge0$, and any relation of the form $i_n=i_{n+k}$ implies $i_n=i_{n+\ell k}$ for all $\ell\ge0$. Similar relations hold with $i$ replaced by $j$ on one or both sides of the equations. A direct case analysis yields the following disjoint possibilities. In each case, the elements $i_1$, $j_1$ are specially marked and the arrows indicate the action of~$\tau$.
\begin{enumerate}
\item[(a)] All relations result from two relations $i_\ell=i_{k+1}$ and $j_n=j_{m+1}$ with $1\le\ell\le k$ and $1\le n\le m$:
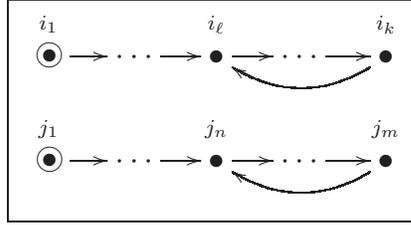
$$\fbox{\ $\xymatrix@R-25pt@C-10pt{
\scriptstyle i_1 && \scriptstyle i_\ell && \scriptstyle i_k \\ 
\circbullet \ar[r]&\cdots\ar[r]&\bullet\ar[r]&\cdots
\ar[r]&\bullet\ar@/^12pt/[ll]  \\
{\vphantom{X}}&&&&\quad \\
\scriptstyle j_1 && \scriptstyle j_n && \scriptstyle j_m \\ 
\circbullet \ar[r]&\cdots\ar[r]&\bullet\ar[r]&\cdots
\ar[r]&\bullet\ar@/^12pt/[ll] \\
{\vphantom{X}}&&&&\\
}$}\\$$
\item[(b)] All relations result from two relations $i_\ell=j_{m+1}$ and $j_n=i_{k+1}$ with $1\le\ell\le k$ and $1\le n\le m$:
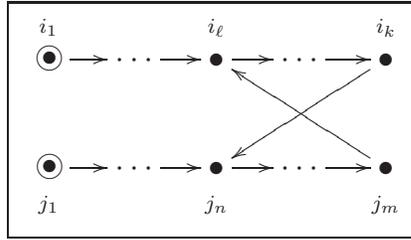
$$\fbox{\ $\xymatrix@R-25pt@C-10pt{
\scriptstyle i_1 && \scriptstyle i_\ell && \scriptstyle i_k \\ 
\circbullet \ar[r]&\cdots\ar[r]&\bullet\ar[r]&\cdots
\ar[r]&\bullet\ar[llddd] \\
{\vphantom{X}}&&&&\\
{\vphantom{X}}&&&&\quad \\
\circbullet \ar[r]&\cdots\ar[r]&\bullet\ar[r]&\cdots
\ar[r]&\bullet\ar[lluuu] \\
\vphantom{j}\mysmash{2pt}{j_1} && \mysmash{2pt}{j_n} && \mysmash{2pt}{j_m} \\ 
}$}\\$$
\item[(c)] 
All relations result from two relations $i_\ell=j_n$ and $i_k=i_{m+1}$ with $1\le\ell<k\le m$ and $1\le n$ and $(\ell,n)\not=(1,1)$:
$$\fbox{\ $
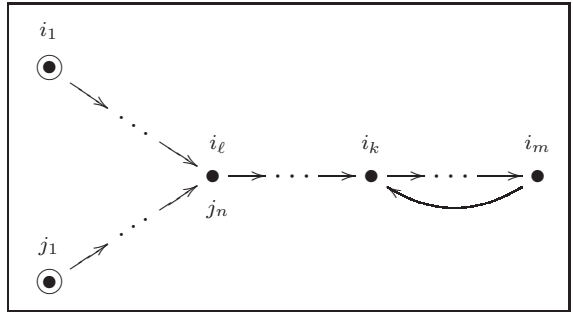
\xymatrix@R-21pt@C-10pt{
\vphantom{x}\mysmash{-1pt}{i_1} &&&&&& \\
\circbullet \ar[dr] &&&&&& \\
&\raisebox{-1.5pt}{$\ddots$}\ar[dr]&\hbox to 0pt{\hss\mysmash{-6pt}{\ i_\ell}\hss}
&&\hbox to 0pt{\hss\mysmash{-6pt}{i_k}\hss}
&&\hbox to 0pt{\hss\mysmash{-6pt}{i_m}\hss}\\
&&\bullet\ar[r]&\cdots\ar[r]&
\bullet\ar[r]&\cdots\ar[r]&\bullet\ar@/^12pt/[ll]\\
\mysmash{-6pt}{j_1} &\raisebox{-1.5pt}{\reflectbox{$\ddots$}}\ar[ur]&
\hbox to 0pt{\hss\mysmash{8pt}{\ j_n}\hss}&&&\\
\circbullet \ar[ur]&&&&&& \quad \\
}$}\qquad\\$$
\end{enumerate}
\end{Class}

\begin{Rem}\label{AllOccur}
\rm By explicit calculation the case (a) with $(\ell,k)=(n,m)=(1,2)$ does not occur as the postcritical orbit of a quadratic morphism over any field. I do not know a conceptual explanation for this fact. By contrast all other finite mapping schemes that were tested in Dau \cite{Dau} occur, say over~$\BC$. Is that true for all other cases as well?
\end{Rem}

%%%%%%%%%%%%%%%%%%%%%%%%%%%%%%%%%%%%%%%%%%%%%%%%%%%%%%%%%%%%%%%%%%%%%%%%%%%%%%%%%%%%%%%%%%%%%

\section{Marked quadratic morphisms}
\label{MarkedQM}

Now we fix a finite mapping scheme $(\Gamma,\tau,i_1,j_1)$. In order to ensure stability we assume that $|\Gamma|\ge3$. (For the remaining cases see Proposition \ref{Unstable}.)

\begin{Def}\label{GammaMarkDef}
A \emph{$\Gamma$-marked quadratic morphism} over $S$ is a quadratic morphism $(C,f,P,Q)$ over $S$ together with a section $s(\gamma)\in C(S)$ for every $\gamma\in\Gamma$ satisfying
\begin{enumerate}
\item[(a)] $s(i_1)=f(P)$ and $s(j_1)=f(Q)$, 
\item[(b)] $f(s(\gamma))=s(\tau(\gamma))$ for any $\gamma\in\Gamma$, and
\item[(c)] $s(\gamma)$, $s(\gamma')$ are fiberwise distinct for any distinct $\gamma,\gamma'\in\Gamma$.
\end{enumerate}
\end{Def}

These conditions guarantee that the map $\gamma\mapsto s(\gamma)$ induces an isomorphism from $\Gamma$ to the postcritical orbit of $f$ in every fiber and give a scheme theoretic version thereof.

\begin{Prop}\label{GammaModuli}
There is a fine moduli scheme $M_\Gamma$ of $\Gamma$-marked quadratic morphisms, which is isomorphic to the joint zero locus of two polynomials in a suitable Zariski open subscheme of $\BA^2\times\Spec\BZ[\half]$.
\end{Prop}

\begin{Proof}
Since $|\Gamma|\ge3$, at least one of $i_2,j_2$ is distinct from both $i_1,j_1$. By symmetry we may assume without loss of generality that $j_2\not=i_1,j_1$. Then for any $\Gamma$-marked quadratic morphism $f$ the section $f^2(Q)$ is disjoint from $f(P)$ and $f(Q)$, and so by (\ref{RamLift}) the section $f(Q)$ is disjoint from $P$ and~$Q$. By the proof of Proposition \ref{StableModuli}, the quadratic morphisms with this property are up to unique isomorphism precisely the $(\BP^1\times S, x\mapsto\frac{x^2+a}{x^2+b},0,\infty)$ for $(a,b)\in (\BA^2\setminus\diag(\BA^1))(S)$.

In the case (a) of the classification \ref{AbsPostCritClass}, there are two relations of the form $i_\ell=i_{k+1}$ and $j_n=j_{m+1}$ which generate all other relations within the mapping scheme~$\Gamma$.
Thus any $\Gamma$-marked quadratic morphism must satisfy the corresponding two relations $f^\ell(P)=f^{k+1}(P)$ and $f^n(Q)=f^{m+1}(Q)$. Conversely, these two relations guarantee that one can attach unique sections $s(\gamma)$ to all $\gamma\in\Gamma$ which satisfy \ref{GammaMarkDef} (a) and (b). The condition \ref{GammaMarkDef} (c) then amounts to finitely many inequalities of the form $f^{n'}(P)\not=f^{m'}(P)$ or $f^{n'}(P)\not=f^{m'}(Q)$
or $f^{n'}(Q)\not=f^{m'}(Q)$ in every fiber. Together this shows that the subscheme of $(\BA^2\setminus\diag(\BA^1))\times\Spec\BZ[\half]$ defined by the two closed conditions $f^\ell(P)=f^{k+1}(P)$ and $f^n(Q)=f^{m+1}(Q)$ and finitely many open conditions is a fine moduli scheme of $\Gamma$-marked quadratic morphisms. The same argument applies to the other cases of the classification~\ref{AbsPostCritClass}.

It remains to see that each of the two closed conditions is represented by a polynomial equation in $(a,b)$. For this observe that in projective coordinates the morphism $f$ is given by $(x:y) \mapsto (x^2+ay^2:x^2+by^2)$, which is well-defined because $a\not=b$ everywhere. Thus the iterate $f^n$ is the well-defined morphism $(x:y) \mapsto (g_n:h_n)$ for certain polynomials $g_n$, $h_n\in \BZ[x,y,a,b]$. Each closed condition is obtained by substituting $(x:y)=(0:1)$ or $(1:0)$ and equating, and therefore means that two well-defined points $(g':h')$ and $(g'':h'')$ in $\BP^1$ for certain polynomials $g'$, $h'$, $g''$, $h''\in \BZ[a,b]$ are equal. But this is equivalent to $g'h''-g''h'=0$, which is a polynomial equation in $a$, $b$ with coefficients in~$\BZ$, as desired.
\end{Proof}

\begin{Thm}\label{GammaModuliQuasiFinite}
The moduli scheme $M_\Gamma$ is quasi-finite over $\Spec\BZ[\half]$, in other words, for any field $k$ the product $M_\Gamma\times\Spec k$ is finite over~$k$.
\end{Thm}

This is the main technical result of this paper.
In characteristic zero it is known as a consequence of Thurston rigidity, which implies that postcritically finite quadratic morphisms cannot be nontrivially deformed. See Douady-Hubbard \cite{DouadyHubbard1993} and Brezin et al. \cite[Thm. 3.6]{BBLPP2000}.
% ; compare also \cite{Ingram2011a}.
% Also, Ingram \cite{Ingram2011a} has a finiteness result in characteristic zero for postcritical polynomials (in one variable) of arbitrary degree, using Silverman's height.
% Baker - deMarco 2012 Prop.2.7: Only countably many PCF polynomials over $\BC$.
After finishing the work on this paper, we learned that this was already proved using different methods by Benedetto-Ingram-Jones-Levy \cite[Cor.~6.3]{BIJL2012}, see also Levy \cite{Levy2012}.

\medskip
Theorem \ref{GammaModuliQuasiFinite} will be proved in Sections \ref{StableCurvesQuad} through~\ref{ProofOfTheorem}. First we deduce some consequences.

\begin{Cor}\label{FieldOfDefinition}
Any $\Gamma$-marked quadratic morphism over a field $k$ can be defined over a subfield which is either finite or a number field.
\end{Cor}

\begin{Proof}
Any $\Gamma$-marked quadratic morphism over $k$ determines a morphism $\Spec k\to M_\Gamma$, and the residue field at the image point on $M_\Gamma$ is a subfield with the desired property.
\end{Proof}

\begin{Cor}\label{GammaModuliFlat}
The moduli scheme $M_\Gamma$ is flat over $\Spec\BZ[\half]$.
\end{Cor}

\begin{Proof}
By Proposition \ref{GammaModuli} and Theorem \ref{GammaModuliQuasiFinite} the moduli scheme $M_\Gamma$ is a complete intersection of codimension $2$ in a regular scheme of finite type. It is therefore Cohen-Macaulay by Matsumura \cite[Thm.~17.3, Thm.~17.4]{Matsumura}.
Since again by Theorem \ref{GammaModuliQuasiFinite} its fibers over $\Spec\BZ[\half]$ are finite, the local flatness criterion \cite[Thm.~23.1]{Matsumura} ensures that $M_\Gamma$ is flat over $\Spec\BZ[\half]$.
\end{Proof}

\begin{Cor}\label{GammaModuliLift}
Any $\Gamma$-marked quadratic morphism $(C,f,P,Q,s)$ over a finite field $k$ of characteristic $p$ can be lifted to characteristic zero. More specifically, there exist a discrete valuation ring $R$ which is finitely generated over $\BZ_{(p)}$, whose residue field $k'$ is a finite extension of $k$ and whose quotient field has characteristic zero, and a $\Gamma$-marked quadratic morphism over $\Spec R$, whose closed fiber is isomorphic to $(C,f,P,Q,s)\times_{\Spec k}\Spec k'$.
\end{Cor}

\begin{Proof}
Let $X$ denote the normalization of the reduced closed subscheme of~$M_\Gamma$. Then the morphism $X\to M_\Gamma$ is finite and surjective; hence the classifying morphism $\Spec k\to M_\Gamma$ lifts to a morphism $\Spec k'\to X$ for some finite extension $k'$ of~$k$. The lift corresponds to a ring homomorphism $R_1\to k'$ for a suitable open affine chart $\Spec R_1\subset X$. This in turn factors through the localization $R_2$ of $R_1$ at some maximal ideal. 
By construction $R_1$ is a normal integral domain, which by Proposition \ref{GammaModuli} and Theorem \ref{GammaModuliQuasiFinite} is finitely generated and flat over~$\BZ$. Thus $R_2$ is a discrete valuation ring that is finitely generated and flat over~$\BZ_{(p)}$.
Let $R$ be an unramified extension of $R_2$ with residue field~$k'$. Then the desired assertion follows by pulling back the universal family over $M_\Gamma$ under the commutative diagram
$$\xymatrix{
\Spec k' \ar@{^{ (}->}[r] \ar[d] & \Spec R \ar[d] \\
\Spec k \ar[r] & M_\Gamma. \\}$$
\vskip-21pt
\end{Proof}

\begin{Rem}\label{RamifiedLift}
\rm In general the ring $R$ in Corollary \ref{GammaModuliLift} may be ramified over $\BZ_{(p)}$. It therefore seems impractical to find such a lift by deformation theory. For a concrete example from Dau \cite{Dau} let $\Gamma$ be the mapping scheme
$$\fbox{$
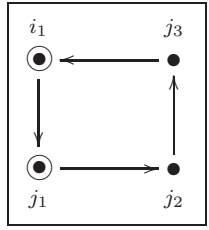
\xymatrix@R-25pt@C-10pt{
\scriptstyle i_1 && \scriptstyle j_3 \\ 
\circbullet\ar[ddd]&&\bullet\ar[ll] \\
{\phantom{X}}&&\\
{\phantom{X}}&& \quad\\
\circbullet\ar[rr]&&\bullet\ar[uuu] \\
\vphantom{x}\mysmash{2pt}{j_1} && \mysmash{2pt}{j_2} \\ 
}$}\qquad$$
which is the case \ref{AbsPostCritClass} (b) with $(\ell,k,n,m)=(1,1,1,3)$. A direct calculation shows that in this case $M_\Gamma \cong \Spec\BZ[\half,a]/(a^2+3a+1)$ with the universal family $(\BP^1\times M_\Gamma, {x\mapsto \frac{x^2+a}{x^2}},0,\infty)$. Here $\BZ[\half,a]/(a^2+3a+1) \cong \BZ[\half,\frac{-3+\sqrt{5}}{2}] = \CO_K[\half]$ where $\CO_K$ is the ring of integers in the quadratic number field $K := \BQ(\sqrt{5})$, which is ramified at the prime $p=5$. 
\end{Rem}

\begin{Rem}\label{NotReduced}
\rm In the same cases as in Remark \ref{RamifiedLift}, the fiber $M_\Gamma\times\Spec\BF_p$ is not reduced. This means that in positive characteristic $\Gamma$-marked quadratic morphisms may possess non-trivial infinitesimal deformations. This makes it impractical to try to prove Theorem \ref{GammaModuliQuasiFinite} by local deformation theory and suggests that the proof requires at least one global argument.
\end{Rem}

\begin{Rem}\label{NotAllOccur}
\rm In general the moduli scheme $M_\Gamma$ is not finite over $\Spec\BZ[\half]$. Indeed, if it were finite over $\Spec\BZ[\half]$, then any point on $M_\Gamma$ in characteristic zero would possess a reduction in any characteristic $\not=2$. But there are cases where $M_\Gamma$ possesses a point in characteristic zero but none in some characteristic $p\not=2$. For a concrete example from Dau \cite{Dau} let $\Gamma$ be the mapping scheme
$$\fbox{\ $
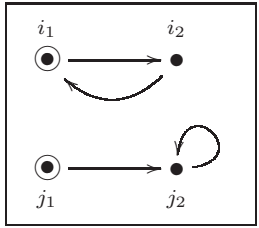
\xymatrix@R-25pt@C-10pt{
\scriptstyle i_1 && \scriptstyle i_2 &\\ 
\circbullet\ar[rr]&&\bullet\ar@/^15pt/[ll] & \\
{\phantom{X}}&&&\\
{\phantom{X}}&&&\\
\circbullet\ar[rr]&&\bullet \ar@(r,u)[] & \\
\vphantom{x}\mysmash{2pt}{j_1} && \mysmash{2pt}{j_2} & \\ 
}$}\qquad$$
which is the case \ref{AbsPostCritClass} (a) with 
$(\ell,k,n,m)=(1,2,2,2)$. A direct calculation shows that in this case $M_\Gamma \cong \Spec\BZ[{\frac{1}{6}}]$ with the universal family $(\BP^1\times M_\Gamma, {x\mapsto \frac{x^2-4}{x^2+2}},0,\infty)$. Here the prime $p=3$ must be excluded in order to have $-4\not\equiv2$. Thus for this~$\Gamma$, there exist $\Gamma$-marked quadratic morphisms in every characteristic $\not=2$, $3$, but none in characteristic~$3$.
\end{Rem}

\begin{Rem}\label{BadReduction}
\rm In the same cases as in Remark \ref{NotAllOccur}, and possibly others, a $\Gamma$-marked quadratic morphism over a number field does not reduce to a $\Gamma$-marked quadratic morphism in characteristic~$p$. In other words, the quadratic morphism has \emph{bad reduction}. In Sections \ref{StableCurvesQuad} and \ref{Exclude} we analyze such situations systematically over arbitrary discrete valuation rings. Although some types of bad reduction may occur, we will prove that a certain kind cannot occur. This will suffice to deduce that any $\Gamma$-marked quadratic morphism over an irreducible reduced curve over a field is constant, and thereby prove Theorem \ref{GammaModuliQuasiFinite}.
\end{Rem}

%%%%%%%%%%%%%%%%%%%%%%%%%%%%%%%%%%%%%%%%%%%%%%%%%%%%%%%%%%%%%%%%%%%%%%%%%%%%%%%%%%%%%%%%%%%%%

\section{Monodromy groups}
\label{MonodromyGroups}

Let $(C,f,P,Q,s)$ be a $\Gamma$-marked quadratic morphism over a connected scheme~$S$. Let $C^\circ$ denote the open subscheme of $C$ obtained by removing the images of the sections $s(\gamma)\colon S\to C$ for all $\gamma\in\Gamma$. For any integer $n\ge0$ set $C_n := C$ as a scheme over~$C$ via the morphism $f^n\colon C_n\to C$, and let $C_n^\circ\subset C_n$ be the open subscheme defined as the fiber product $C_n\times_CC^\circ$. Then by construction $C_n^\circ\to C^\circ$ is a finite \'etale covering of degree~$2^n$. The morphism $f$ induces transition morphisms $C_n^\circ\to C_{n-1}^\circ$ for all $n>0$ and hence an inverse system
$$\ldots\onto C_2^\circ\onto C_1^\circ\onto C_0^\circ=C^\circ.$$
% Its inverse limit $C_\infty^\circ$ is then a profinite \'etale covering of~$C^\circ$.
Let $\bar c_0$ be any geometric point of~$C^\circ$. For any $n\ge0$ let $T_n$ denote the set of $2^n$ geometric points of $C_n^\circ$ above~$\bar c_0$. Let $T$ denote the infinite directed graph with set of vertices $\coprod_{n\ge0}T_n$, where any vertex $\bar c_n\in T_n$ for $n>0$ is connected by an edge towards $f(\bar c_n)\in T_{n-1}$. By construction this is a regular rooted binary tree.
Since by assumption $S$ is connected, so is $C_n^\circ$ for all $n$, and so the inverse system is determined up to isomorphism by the monodromy representation of the \'etale fundamental group on~$T$, that is, by the natural homomorphism 
\UseTheoremCounterForNextEquation
\begin{equation}\label{MonoRep}
\pi_{1,\et}(C^\circ,\bar c_0) \to \Aut(T).
\end{equation}
Suppose that $\bar c_0$ lies above the geometric point $\bar s$ of~$S$. Then we are interested in the restriction of the above homomorphism to the geometric \'etale fundamental group 
\UseTheoremCounterForNextEquation
\begin{equation}\label{RelMonoRep}
\rho\colon \pi_{1,\et}(C^\circ_{\bar s},\bar c_0) \to \Aut(T).
\end{equation}

\begin{Prop}\label{SameRelMono}
Let $\bar c_0^{\,\prime}$ be a geometric point of $C_0$ lying above another geometric point $\bar s'$ of~$S'$. Let $T'$ be the regular rooted binary tree constructed with $\bar c_0^{\,\prime}$ in place of $\bar c_0$, and let 
$$\rho'\colon \pi_{1,\et}(C^\circ_{\bar s'},\bar c_0^{\,\prime}) \to \Aut(T')$$
be the associated monodromy representation. Then there is an isomorphism of rooted trees $T\cong T'$ which induces an isomorphism $\image(\rho)\cong\image(\rho')$.
\end{Prop}

\begin{Proof}
The target $\Aut(T)$ and hence the images of $\rho$ and $\rho'$ are pro-$p$-groups for $p=2$. Since $S$ is a scheme over $\Spec\BZ[\half]$, the assertion is a direct consequence of Grothendieck's specialization theorem for the tame fundamental group, see SGA1 \cite[exp.\ XIII, Cor.\ 2.12]{SGA1}.
\end{Proof}

\begin{Cor}\label{RelMonoChar0}
Let $\bar G\subset\Aut(T)$ be the image of the representation (\ref{RelMonoRep}) of the geometric monodromy group for a $\Gamma$-marked quadratic morphism over any field. Then up to an isomorphism of rooted trees $\bar G$ is equal to the image of the geometric monodromy representation for some $\Gamma$-marked quadratic morphism over~$\BC$.
\end{Cor}

\begin{Proof}
Direct consequence of Proposition \ref{SameRelMono} and Corollaries \ref{FieldOfDefinition} and \ref{GammaModuliLift}.
\end{Proof}

Finally suppose that $\bar s$ is a $\BC$-valued point. Then the \'etale fundamental group $\pi_{1,\et}(C^\circ_{\bar s},\bar c_0)$ is naturally isomorphic to the profinite completion of the topological fundamental group $\pi_1(C^\circ(\BC),\bar c_0)$. Thus the image of the former is simply the closure in the profinite topology of $\Aut(T)$ of the image of $\pi_1(C^\circ(\BC),\bar c_0)$. But the latter image can be described by explicit generators as a self-similar subgroup of $\Aut(T)$ as, say, in 
Bartholdi-Nekrashevych \cite{Bartholdi-Nekrashevych-2008},
Dau \cite{Dau}, 
Grigorchuk et al.~\cite{Grigorchuk-et-al-2007},
Nekrashevych \cite[Ch.5]{Nekrashevych-2005}, \cite{Nekrashevych-2009}.

%%%%%%%%%%%%%%%%%%%%%%%%%%%%%%%%%%%%%%%%%%%%%%%%%%%%%%%%%%%%%%%%%%%%%%%%%%%%%%%%%%%%%%%%%%%%%

% \newpage
\Bigskip
\noindent{\bf\Large Part II: Finiteness}
\label{Finite}
\bigskip

The rest of this paper is devoted to proving Theorem \ref{GammaModuliQuasiFinite}, which seems surprisingly difficult. 

%%%%%%%%%%%%%%%%%%%%%%%%%%%%%%%%%%%%%%%%%%%%%%%%%%%%%%%%%%%%%%%%%%%%%%%%%%%%%%%%%%%%%%%%%%%%%

\section{Stable marked curves of genus zero}
\label{StableCurves}

We begin by briefly reviewing the material on stable marked curves that we will need below; see Deligne-Mumford \cite{DeligneMumford}, Knudsen, \cite{Knudsen1983a}, Gerritzen et al.~\cite{GHP}, Keel \cite{Keel1992}.

Let $R$ be a discrete valuation ring with quotient field~$K$, uniformizer~$\pi$, and residue field $k=R/R\pi$. Let $C$ be a smooth connected curve of genus zero over~$K$. Let $I$ be a finite set of cardinality $\ge3$, and consider an injective map $s\colon I\into C(K)$, $i\mapsto s(i)$. Then $(C,s)$ is a \emph{smooth marked curve of genus zero} over~$K$.
By the general theory of stable marked curves it possesses a natural extension to a \emph{stable marked curve $(X,s)$ of genus zero} over~$\Spec R$, as follows. 

First $X\to\Spec R$ is a projective and flat morphism with generic fiber~$C$. Let $X^\sm\subset X$ denote the open locus where the morphism is smooth, and let $X^\sm_0\subset X_0$ denote the corresponding closed fibers. Then the complement $X\setminus X^\sm$ is a finite subset of~$X_0$, and at each point in it $X$ is \'etale locally isomorphic to $\Spec R[x,y]/(xy-\pi^n)$ over $\Spec R$ for some integer $n\ge1$. 
Next the markings $s(i)\in C(K)$ extend to pairwise disjoint sections $s(i)\colon \Spec R\to X^\sm$ and hence induce an injection $s_0\colon I\into X^\sm_0(k)$.
Finally, the data $(X,s)$ is globally stable in the sense that for every irreducible component $Y$ of~$X_0$, we have 
\UseTheoremCounterForNextEquation
\begin{equation}\label{StabCond1}
\bigl|\{i\in I\mid s_0(i) \in Y(k)\}\bigr| + \bigl|Y\setminus X_0^\sm\bigr| \ge3.
\end{equation}
The extension $(X,s)$ thus characterized is unique up to unique isomorphism.
% Note: The conditions imply that every irreducible component and every singular point of $X_0$ is defined over~$k$.

The combinatorial structure of the special fiber $X_0$ is described in terms of its \emph{dual tree}. This is the finite graph $T$ whose set of vertices $V(T)$ is in bijection with the set of irreducible components of $X_0$ and whose set of edges with the set of singular points of~$X_0$, where any singular point corresponds to an edge between the two irreducible components in which it is contained. In our case the fact that $C$ has genus zero implies that this graph is actually a tree. In particular, it possesses no multiple edges and no edge connecting a vertex to itself. Also, every irreducible component of $X_0$ is smooth of genus zero over~$k$.
In addition, the marking $s_0\colon I\to X^\sm_0(k)$ induces a map 
\UseTheoremCounterForNextEquation
\begin{equation}\label{MarkOnTree}
I\to V(T)
\end{equation}
sending $i$ to the unique irreducible component $Y$ of $X_0$ with $s_0(i) \in Y(k)$. We view this map as a marking on the tree~$T$, and by abuse of notation we denote it again by~$s$. Of course this map is no longer injective, and the stability condition (\ref{StabCond1}) translates into the following condition for every vertex $t\in V(T)$:
\UseTheoremCounterForNextEquation
\begin{equation}\label{StabCond2}
\bigl|s^{-1}(t)\bigr| + \bigl|\{\hbox{edges at $t$}\}\bigr| \ge3.
\end{equation}

%%%%%%%%%%%%%%%%%%%%%%%%%%%%%%%%%%%%%%%%%%%%%%

\medskip
Now consider a subset $I'\subset I$ which is also of cardinality $\ge3$. Then we can apply everything above with the marking $s|_{I'}$ in place of~$s$, obtaining another stable marked curve $(X',s')$ over~$\Spec R$. Its relation with $(X,s)$ is described as follows.

First, the identity morphism on $C$ extends to a unique morphism $\kappa\colon X\onto X'$ satisfying $\kappa\circ s|_{I'} = s'$. This morphism contracts the irreducible components of $X_0$ which become unstable, that is, which violate condition (\ref{StabCond1}) when $s$ is replaced by $s|_{I'}$. After these unstable irreducible components have been contracted, the number of singular points of the special fiber may have decreased, and so other irreducible components may have become unstable as well. After finitely many steps, however, the result is stable and is the $X'$ obtained from the marking $s|_{I'}$. 

Next let $T'$ be the dual tree of the special fiber $X_0'$ of~$X'$. Mapping each irreducible component of $X_0'$ to its proper transform in $X_0$ defines a natural injection $V(T')\into V(T)$. We identify $V(T')$ with its image and call it the set of vertices which \emph{survive in $T'$}. 
It must be noted that the resulting diagram
\UseTheoremCounterForNextEquation
\begin{equation}\label{StabilizationDiag}
\vcenter{\xymatrix{
V(T') \ar@{^{ (}->}[r] & V(T) \\
I' \ar@{^{ (}->}[r] \ar[u]^{s'} & I  \ar[u]^{s}\\}}
\end{equation}
does not commute in general. Instead, the map $s'$ in it is obtained from $s$ by Proposition \ref{StabilizationProp} (c) below.

When we say that a vertex $t''$ of $T$ lies \emph{between} two vertices $t$ and~$t'$, we always mean that $t''$ lies on the shortest path from $t$ to~$t'$, including $t$ and~$t'$, and even allowing $t''=t=t'$. When we say that $t''$ lies \emph{strictly between} $t$ and~$t'$, we require in addition that $t''\not=t,t'$.

\begin{Prop}\label{StabilizationProp}
\begin{enumerate}
\item[(a)] A vertex $t\in V(T)$ survives in $T'$ if and only if 
% $$\bigl|s^{-1}(t)\cap I'\bigr| + \bigl|\{\hbox{$S$ connected component of $T\setminus\{t\}$ with $s^{-1}(V(S))\cap I'\not=\emptyset$}\}\bigr| \ge3.$$
$$\bigl|s^{-1}(t)\cap I'\bigr| + \left|\left\{
\begin{array}{c}
\hbox{$S$ connected component of $T\setminus\{t\}$} \\
\hbox{such that $s^{-1}(V(S))\cap I'\not=\emptyset$}
\end{array}
\right\}\right| \ge3.$$
\item[(b)] Two distinct vertices $t$, $t'\in V(T')$ are connected by an edge in~$T'$ if and only if 
every vertex strictly between $t$ and $t'$ does not survive in~$T'$.
\item[(c)] For any $i\in I'$ the $s'(i)$ is the unique vertex in $V(T')$ whose distance to $s(i)$ in $T$ is minimal.
In particular we have $s'(i)=s(i)$ if $s(i)\in V(T')$.
\end{enumerate}
\end{Prop}

\begin{Proof}
Induction on the cardinality of $I\setminus I'$. If $I'=I$, then for any $t\in V(T)$ any connected component $S$ of $T\setminus\{t\}$ contains an end vertex of $T$ and hence a marked point. Thus the condition in (a) is satisfied by the stability of $(T,s)$, while the assertions (b) and (c) hold trivially in this case.

If $I'\not=I$ set $I'' := I\setminus\{i_0\}$ for some $i_0\in I\setminus I'$, let $(X'',s'')$ be the stable marked curve obtained by stabilizing $(X,s|_{I''})$, and let $(T'',s'')$ be the marked dual tree of its special fiber~$X_0''$. Then by the induction hypothesis the assertions (a) through (c) hold for obtaining $(T',s')$ from $(T'',s'')$. 
If all vertices of $T$ survive in~$T''$, the assertions for obtaining $(T',s')$ from $(T,s)$ follow directly. 

Otherwise $T''$ is obtained from $T$ by dropping the vertex $t_0 := s(i_0)$, and the corresponding irreducible component $Y$ of $X_0$ is contracted to a point. By the stability of $X_0$ this can happen only in one of two cases: 

Either $t_0$ is an end vertex, i.e., a leaf, of~$T$ and is marked by precisely one other point $i_1\in I''$. Then $Y$ is contracted to a smooth point of~$X''_0$ which becomes marked by $s''(i_1)$. This means that the marking $s(i_1)$ of $T$ is moved to the unique neighboring vertex in~$T$. One easily checks that the assertions (a) through (c) for obtaining $(T',s')$ from $(T,s)$ now follow directly from those for obtaining $(T',s')$ from $(T'',s'')$. 

Or $t_0$ has precisely two neighbors in $T$ and no other marking in~$I''$. Then the two neighboring vertices correspond to the unique two other irreducible components of $X_0$ which meet~$Y$. These are disjoint in~$X_0$, but their images in $X_0''$ meet at the point obtained by contracting~$Y$. Thus as the vertex $t_0$ is dropped in~$T''$, the two edges connecting it to its neighbors in $T$ are replaced by a single edge in $T''$ connecting these neighbors directly. Again one easily checks that the assertions (a) through (c) for obtaining $(T',s')$ from $(T,s)$ follow directly from those for obtaining $(T',s')$ from $(T'',s'')$. 
\end{Proof}

%%%%%%%%%%%%%%%%%%%%%%%%%%%%%%%%%%%%%%%%%%%%%%%%%%%%%%%%%%%%%%%%%%%%%%%%%%%%%%%%%%%%%%%%%%%%%

\section{Stable marked curves with an automorphism of order two}
\label{StableAut}

Keeping the notations of the preceding section, we now assume that the characteristic of $k$ and hence of $K$ is different from~$2$. Let $\sigma$ be an automorphism of order $2$ of $C$ over~$K$. Then $\sigma$ possesses precisely two fixed points over an algebraic closure of~$K$. We assume that both of these are marked points, i.e., that they are equal to $s(i_0)$ and $s(j_0)$ for certain distinct elements $i_0$, $j_0\in I$. 
Let $\sigma$ also denote a permutation of order $2$ of $I$, and assume that $\sigma\mycirc s = s\mycirc\sigma$. Then the injectivity and $\sigma$-equivariance of $s$ implies that $i_0$, $j_0$ are precisely the fixed points of $\sigma$ on~$I$. 
The uniqueness of the stable extension $(X,s)$ implies that $\sigma$ extends to a unique automorphism of order $2$ of $X$ over~$\Spec R$, which we again denote by~$\sigma$. The aim of this section is to describe the action of~$\sigma$ on the closed fiber and to analyze the quotient $X/\gensigma$. This can be done in a relatively explicit way.

\medskip
Without loss of generality we may assume that $C=\BP^1_K$ with $\sigma(x)=-x$ and $s(i_0)=0$ and $s(j_0)=\infty$. Let $\infty>n_1>\ldots>n_r>-\infty$ be the possible orders $\ord_\pi(s(i))$ for all $i\in I\setminus\{i_0,j_0\}$. Since $|I|\ge3$, there is at least one, and so $r\ge1$. Define
$$U_\nu\ :=\ \left\{
\begin{array}{ll}
% \Spec R\bigl[{\textstyle\frac{x}{\pi^{n_1}}}\bigr] 
\Spec R\,[\,x/\pi^{n_1}] 
& \hbox{if\ \ $\nu=0$,} \\[6pt]
% \Spec R\bigl[{\textstyle\frac{x}{\pi^{n_{j+1}}},\frac{\pi^{n_j}}{x}}\bigr]
\Spec R\,[\,x/\pi^{n_{\nu+1}},\pi^{n_\nu}\!/x\,]
& \hbox{if\ \ $0<\nu<r$,} \\[6pt]
% \Spec R\bigl[{\textstyle\frac{\pi^{n_r}}{x}}\bigr]
\Spec R\,[\,\pi^{n_r}\!/x\,]
& \hbox{if\ \ $\nu=r$.}
\end{array}
\right.$$
For any $1\le\nu\le r$ we glue $U_{\nu-1}$ and $U_\nu$ along the common open subscheme 
$$U_{\nu-1}\cap U_\nu = \Spec R\,[\,(x/\pi^{n_\nu})^{\pm1}],$$
obtaining a projective flat curve $Z$ over $\Spec R$. For all $1\le \nu\le r$ let $Y_\nu$ denote the closure in $Z$ of the closed fiber
$$(U_{\nu-1}\cap U_\nu)_0 \ =\ 
% \Spec R\bigl[\bigl({\textstyle\frac{x}{\pi^{n_j}}}\bigr)^{\pm1}\bigr]/(\pi)
\Spec R\,[\,(x/\pi^{n_\nu})^{\pm1}]/(\pi)
\ \cong\ \BP^1_k\setminus\{0,\infty\}$$
of $U_{\nu-1}\cap U_\nu$. Then $Y_\nu\cong\BP^1_k$, and these are precisely the irreducible components of the closed fiber $Z_0$ of~$Z$. They are arranged in sequence, each meeting only the previous and the next one in an ordinary double point. The automorphism $\sigma\colon x\mapsto-x$ extends to an automorphism of each~$U_\nu$, and since $\sigma(x/\pi^{n_\nu}) = - x/\pi^{n_\nu}$, the automorphism induced on $Y_\nu\cong\BP^1_k$ has the form $\xi\mapsto-\xi$. As the residue characteristic is $\not=2$, it follows that $\sigma$ has exactly two fixed points on each~$Y_\nu$. Varying~$\nu$, it follows that the fixed points of $\sigma$ on $Z_0$ are precisely the singular points of $Z_0$ and one additional point on each of $Y_1$ and~$Y_r$. The last two are actually the reductions of the points $s(i_0)=0$ and $s(j_0)=\infty$, which extend to sections of $U_1$ and~$U_r$ that we again denote by $s(i_0)$ and $s(j_0)$, respectively.

For any $i\in I\setminus\{i_0,j_0\}$, by construction there is a unique $1\le \nu\le r$ with $\ord_\pi(s(i))=n_\nu$. Then $s(i)/\pi^{n_\nu}$ lies in~$R^\times$; hence the point $s(i)\in C(K)$ extends to a section $\Spec R\to U_{\nu-1}\cap U_\nu$. In the closed fiber this section meets~$Y_\nu$, and since $s(i)/\pi^{n_\nu} \mod (\pi) \not=0$,~$\infty$, the resulting point of $Y_\nu$ is not fixed by~$\sigma$. In particular the section lands in the smooth locus $Z^\sm$ of~$Z$ and is disjoint from the sections $s(i_0)$ and $s(j_0)$. For simplicity we denote this section again by~$s(i)$. Altogether we have thus extended the given points $s(i)\in C(K)$ to sections $s(i)\colon\Spec R\to Z^\sm$ for all $i\in I$. Since $\sigma\mycirc s = s\mycirc\sigma$ in the generic fiber, the same relation holds for the extended sections as well.

The choice of the $n_\nu$ also implies that for every $1\le \nu\le r$, there exists $i\in I\setminus\{i_0,j_0\}$ with $\ord_\pi(s(i))=n_\nu$. The corresponding section $s(i)$ in the special fiber then lands in the irreducible component~$Y_\nu$. Since the resulting point of $Y_\nu$ is not fixed by~$\sigma$, its image under $\sigma$ is a different marked point in~$Y_\nu$. In~$Y_1$, respectively~$Y_r$, there is also the third distinct marked point $s(i_0)$, respectively $s(j_0)$, and the remaining irreducible components contain two singular points of~$Z_0$. Thus the closed fiber $Z_0$ with the marking $s$ satisfies the stability condition (\ref{StabCond1}).

All this shows that $(Z,s)$ is a stable marked curve of genus zero, \emph{except} that some (and possibly very many) sections $s(i)$ and $s(i')$ for $i$, $i'\in I\setminus\{i_0,j_0\}$ with $i'\not\in\{i,\sigma(i)\}$ may meet in the special fiber. Let $S$ be the finite set of points in $Z_0^\sm$ where this happens. Then by the general theory of stable curves, the stable extension $(X,s)$ is obtained from $(Z,s)$ by blowing up some ideal centered in $S$ while leaving the rest of $Z$ unchanged. 
The proper transforms $\tilde Y_\nu$ of $Y_\nu$ in~$X$ are then still arranged in the same way as the~$Y_\nu$.
But the sections $s(i)$ reducing to points in $S$ are moved to new irreducible components in the exceptional fibers of the blowup. By the uniqueness of $X$ everything is still invariant under the action of~$\sigma$. Also, each new irreducible component is disjoint from its $\sigma$-conjugate, because its image in $Z_0$ is a point in $S$ which is not fixed by~$\sigma$. 

The following sketch shows what may typically happen. The irreducible fibers of the closed fiber are indicated in solid lines, the marked sections in dotted lines:
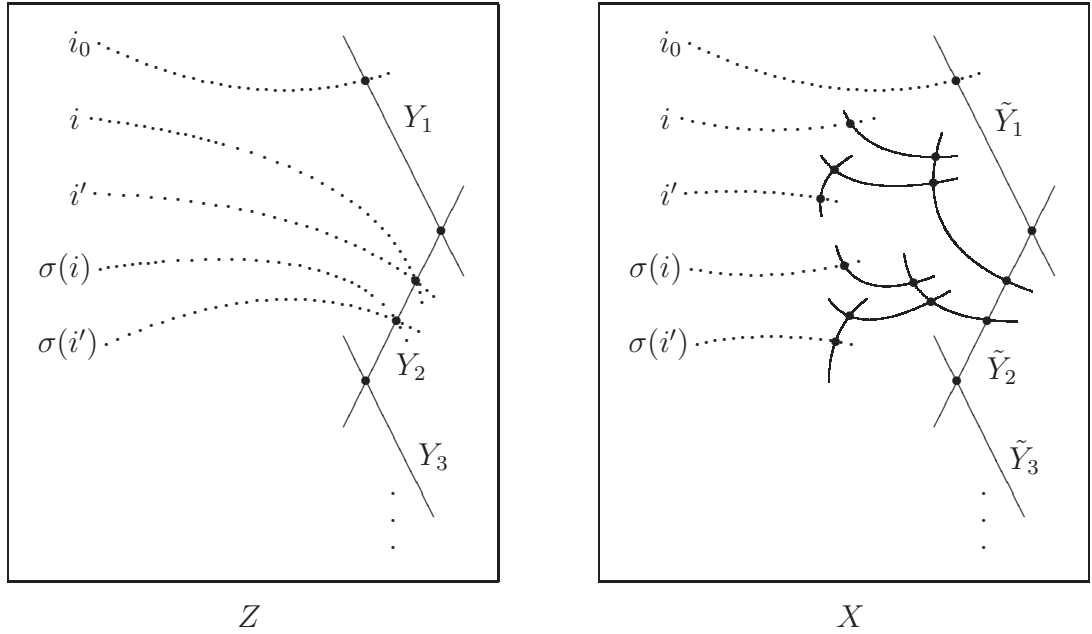
\begin{center}
\smallskip
\def\myotherbullet{\scriptstyle\bullet}
\parbox[t]{.4\textwidth}
{\setlength{\fboxsep}{10pt}
\fbox{\xy 
(-4,40)*{}; 
% (-7,-11)*{}; 
(54,100)*{}; 
(39.9,100)*{\myotherbullet};
(53,74)*{};(37,106)*{};**\crv{};
(47,95)*{Y_1};
(50,80)*{\myotherbullet};
(53,86)*{};(37,54)*{};**\crv{};
(46,62)*{Y_2};
% (40.1,0.05)*{\myotherbullet};
(49,42)*{};(37,66)*{};**\crv{};
(49,50)*{Y_3};
(45,37)*{};(45,47)*{};**\crv{~**=<8pt>{.}};
% (50,20)*{};(50,40)*{};**\crv{~**=<10pt>{.}};
% (49,18)*{};(37,-6)*{};**\crv{};
% (49,2)*{Y_r};
(40,60)*{\myotherbullet};
(46.65,73.3)*{\myotherbullet};
(44.05,68.1)*{\myotherbullet};
(4.5,105)*{};(43,101)*{};**\crv{~*=<4pt>{}~**=<4pt>{.}(25,95)};
(2,105)*{i_0};
(3.5,95)*{};(47.5,70.6)*{};**\crv{~*=<4pt>{}~**=<4pt>{.}(43,90)};
(1.3,95)*{i};
(4,85)*{};(49,71.4)*{};**\crv{~*=<4pt>{}~**=<4pt>{.}(35,84)};
(1.8,85)*{i'};
(4.7,75)*{};(45.4,65.5)*{};**\crv{~*=<4pt>{}~**=<4pt>{.}(40,80)};
(0,75)*{\sigma(i)};
(5.5,65)*{};(47,66.7)*{};**\crv{~*=<4pt>{}~**=<4pt>{.}(30,76)};
(0.4,65)*{\sigma(i')};
% (46.6,13.2)*{\myotherbullet};
% (5,20)*{};(50,12)*{};**\crv{~*=<4pt>{}~**=<4pt>{.}(26,22)};
% (2,20)*{i''};
% (43.95,7.8)*{\myotherbullet};
% (6.3,10.5)*{};(47,7)*{};**\crv{~*=<4pt>{}~**=<4pt>{.}(26,14)};
% (0.8,10)*{\sigma(i'')};
% (4.5,-2)*{};(44,-.2)*{};**\crv{~*=<4pt>{}~**=<4pt>{.}(20,2)};
% (2,-3)*{i_2};
\endxy}\\[6pt]
\strut\hskip87pt $Z$}
\qquad\quad
\parbox[t]{.4\textwidth}
{\setlength{\fboxsep}{10pt}
\fbox{\xy 
(-4,40)*{}; 
% (-7,-11)*{}; 
(54,100)*{}; 
(39.9,100)*{\myotherbullet};
(53,74)*{};(37,106)*{};**\crv{};
(47,95)*{\tilde Y_1};
(50,80)*{\myotherbullet};
(53,86)*{};(37,54)*{};**\crv{};
(46,62)*{\tilde Y_2};
% (40.1,0.05)*{\myotherbullet};
(49,42)*{};(37,66)*{};**\crv{};
(49,50)*{\tilde Y_3};
(45,37)*{};(45,47)*{};**\crv{~**=<8pt>{.}};
% (50,20)*{};(50,40)*{};**\crv{~**=<10pt>{.}};
% (49,18)*{};(37,-6)*{};**\crv{};
% (49,2)*{\tilde Y_r};
(40,60)*{\myotherbullet};
(46.65,73.3)*{\myotherbullet};
(50,72)*{};(38,93)*{};**\crv{(33,78)};
(37.2,89.85)*{\myotherbullet};
(25,96)*{};(40,90)*{};**\crv{(28,89)};
(36.9,86.35)*{\myotherbullet};
(22,90)*{};(40,87)*{};**\crv{(26,84)};
(23.75,88.1)*{\myotherbullet};
(26,90)*{};(22,82)*{};**\crv{(21,87)};
(44.05,68.1)*{\myotherbullet};
(48,68)*{};(33,77)*{};**\crv{(34,68)};
(34.25,73.1)*{\myotherbullet};
(24,78)*{};(37,74)*{};**\crv{(26,70)};
(36.6,70.6)*{\myotherbullet};
(23,71)*{};(39,72)*{};**\crv{(27,65)};
(25.75,68.65)*{\myotherbullet};
(28,71)*{};(23,60)*{};**\crv{(23,67)};
(4.5,105)*{};(43,101)*{};**\crv{~*=<4pt>{}~**=<4pt>{.}(25,95)};
(2,105)*{i_0};
(25.9,94.3)*{\myotherbullet};
(3.5,95)*{};(29,95)*{};**\crv{~*=<4pt>{}~**=<4pt>{.}(15,92)};
(1.3,95)*{i};
(21.9,84.3)*{\myotherbullet};
(4,85)*{};(24,84)*{};**\crv{~*=<4pt>{}~**=<4pt>{.}(15,86)};
(1.8,85)*{i'};
(25.1,75.4)*{\myotherbullet};
(4.7,75)*{};(27,76)*{};**\crv{~*=<4pt>{}~**=<4pt>{.}(15,73)};
(0,75)*{\sigma(i)};
(23.9,65.3)*{\myotherbullet};
(5.5,65)*{};(26,65)*{};**\crv{~*=<4pt>{}~**=<4pt>{.}(15,67)};
(0.4,65)*{\sigma(i')};
% (46.6,13.2)*{\myotherbullet};
% (5,20)*{};(50,12)*{};**\crv{~*=<4pt>{}~**=<4pt>{.}(26,22)};
% (2,20)*{i''};
% (43.95,7.8)*{\myotherbullet};
% (6.3,10.5)*{};(47,7)*{};**\crv{~*=<4pt>{}~**=<4pt>{.}(26,14)};
% (0.8,10)*{\sigma(i'')};
% (4.5,-2)*{};(44,-.2)*{};**\crv{~*=<4pt>{}~**=<4pt>{.}(20,2)};
% (2,-3)*{i_2};
\endxy}\\[6pt]
\strut\hskip90pt $X$}\quad\strut
\end{center}

% \bigskip
The following proposition summarizes some of the main information gathered so far:

\begin{Prop}\label{SigmaActProp}
\begin{enumerate}
\item[(a)] The fixed points of $\sigma$ in $X_0$ are precisely the reductions of the sections $s(i_0)$ and $s(j_0)$ and the double points of $X_0$ which separate them.
\item[(b)] Any irreducible component $Y$ of $X_0$ is either equal to $\sigma(Y)$ or disjoint from $\sigma(Y)$.
\item[(c)] An irreducible component $Y$ of $X_0$ is equal to $\sigma(Y)$ if and only if it contains a fixed point of~$\sigma$. The automorphism induced by $\sigma$ on it is then non-trivial.
\end{enumerate}
\end{Prop}

Now consider the quotient $\barX := X/\gensigma$, which exists because $X$ is projective over $\Spec R$. Let $\varpi\colon X\onto\barX$ denote the natural morphism. Set $\bar I :=I/\gensigma$, and for every orbit $\bar i = \{i,\sigma(i)\} \in \bar I$ consider the section $\bar s(\,\bar i\,) := \varpi\circ s(i)\colon\Spec R\to\barX$.

\begin{Prop}\label{QuotientProp}
\begin{enumerate}
\item[(a)] The pair $(\barX,\bar s)$ is a stable marked curve over $\Spec R$.
\item[(b)] For any double point $P$ of $X_0$ which is fixed by~$\sigma$ and where $X$ is \'etale locally isomorphic to $\Spec R[y,z]/(yz-\pi^n)$, the quotient $\barX$ is \'etale locally isomorphic to $\Spec R[u,v]/(uv-\pi^{2n})$ at~$\varpi(P)$, for some $n>0$.
\end{enumerate}
\end{Prop}

\begin{Proof}
Away from the fixed points of $\sigma$ the morphism $\varpi$ is \'etale. Thus away from the images of these fixed points, the quotient $\barX$ has the same mild singularities as~$X$. Moreover, the sections $s(i)$ for all $i\in I\setminus\{i_0,j_0\}$ land in the smooth locus $X^\sm$ and are pairwise disjoint from $s(i_0)$ and $s(j_0)$ and each other. Thus by Proposition \ref{SigmaActProp} (a) they avoid all fixed points of~$\sigma$, and so the corresponding sections $\bar s(\,\bar i\,)$ land in the smooth locus $\barX^\sm$ and are again disjoint from $\bar s(\,\bar i_0\,)$ and $\bar s(\,\bar j_0\,)$ and each other.

Next, by the above construction $X$ is locally isomorphic to $U_1=\Spec R\,[\,x/\pi^{n_1}]$ along the section $s(i_0)$, with the action $\sigma(x/\pi^{n_1}) = - x/\pi^{n_1}$. Thus the substitution $u=x^2/\pi^{2n_1}$ shows that $\barX$ is locally isomorphic to $U_1=\Spec R[u]$ along the section $\bar s(\,\bar i_0\,)$, and $\bar s(\,\bar i_0\,)$ is given by $u=0$. Thus $\bar s(\,\bar i_0\,)$ lands in the smooth locus~$\barX^\sm$. By symmetry the same is true for the section $\bar s(\,\bar j_0\,)$, which remains disjoint from $\bar s(\,\bar i_0\,)$.
Thus all sections $\bar s(\,\bar i\,)$ for $\bar i\in\bar I$ land in $\barX^\sm$ and are pairwise disjoint.

At the remaining fixed points $X$ is by construction locally isomorphic to 
$$U_\nu = \Spec R\,[\,x/\pi^{n_{\nu+1}},\pi^{n_\nu}\!/x\,]$$ 
for $0<\nu<r$. 
Since $\sigma(x)=-x$, it follows that at the image of this fixed point $\barX$ is locally isomorphic to
$$U_\nu/\gensigma = \Spec R\,[\,x^2/\pi^{2n_{\nu+1}},\pi^{2n_\nu}\!/x^2\,].$$
With the substitutions $y=x/\pi^{n_{\nu+1}}$ and $z=\pi^{n_\nu}\!/x$ the first chart 
is isomorphic to $\Spec R\,[\,y,z\,]/(yz-\pi^n)$ where $n:=n_\nu-n_{\nu+1}>0$, and with $u=y^2$ and $v=z^2$ the second chart becomes isomorphic to $\Spec R[u,v]/(uv-\pi^{2n})$. Thus $\barX$ has the required local form, and we have also proved assertion (b).

It remains to verify the stability condition (\ref{StabCond1}) for $(\barX,\bar s)$. By construction the irreducible components of the closed fiber $\barX_0$ of $\barX$ are precisely the images $\varpi(Y)$ of the irreducible components $Y$ of~$X_0$. 

If $Y\not=\sigma(Y)$, then $Y$ is disjoint from $\sigma(Y)$ by Proposition \ref{SigmaActProp} (b), and so $Y$ maps isomorphically to its image in~$\barX$. Moreover, for any section $s(i)$ meeting $Y$ the conjugate $\sigma(s(i))$ meets $\sigma(Y)$ and therefore not~$Y$. Thus the number of sections $s(i)$ meeting $Y$ is equal to the number of sections $\bar s(\,\bar i\,)$ meeting $\varpi(Y)$. This shows that the left hand sides of (\ref{StabCond1}) for $Y$ and $\varpi(Y)$ are equal. The stability at $Y$ thus implies the stability at $\varpi(Y)$.

If $Y=\sigma(Y)$, then $Y$ is one of the irreducible components~$\tilde Y_\nu$. In it, each of the two fixed points of $\sigma$ either arises from the section $s(i_0)$ or $s(j_0)$ or is a singular point of~$X_0$. Thus each of them contributes $1$ to the left hand side of (\ref{StabCond1}) for $(X,s)$ and hence also for $(\barX,\bar s)$. 
By the stability at $Y_\nu$ there must be at least one other contribution for $(X,s)$, that is, another section $s(i)$ with $i\in I\setminus\{i_0,j_0\}$ or another irreducible component not among $Y_1,\ldots,Y_r$ which meets~$Y_\nu$. Taking quotients, this section or irreducible component yields another section or irreducible component which contributes $1$ to the left hand side of (\ref{StabCond1}) for $\varpi(Y)$ in $(\barX,\bar s)$. Again the stability at $\varpi(Y)$ follows.

Altogether the proposition is thus proved. 
(Note: The stability of $(\barX,\bar s)$ implies that $|\bar I|\ge3$, so we did not need to prove this separately. It follows, however, at once from the fact that 
$I$ contains two fixed points of~$\sigma$ and that $|I|\ge3$.)
\end{Proof}

%%%%%%%%%%%%%%%%%%%%%%%%%%%%%%%%%%%%%%%%%%%%%%

\medskip
Continuing the analysis, let now $T$ be the dual tree of $X_0$ endowed with the marking $s\colon I\to V(T)$, as in Section \ref{StableCurves}. By functoriality $\sigma$ induces an automorphism of~$T$, which we (yet again!) denote by~$\sigma$, and which satisfies $\sigma^2=\id$ and $\sigma\mycirc s=s\mycirc\sigma$. Let $t_1,\ldots,t_r$ denote the vertices of $T$ which correspond to the irreducible components $\tilde Y_1,\ldots,\tilde Y_r$ of~$X_0$. 

\begin{Prop}\label{SigmaTreeProp}
\begin{enumerate}
\item[(a)] The fixed points of $\sigma$ on $V(T)$ are precisely the vertices $t_1,\ldots,t_r$.
\item[(b)] The vertices $t_1,\ldots,t_r$ are distinct, connected in the given order by a string of edges, and satisfy $s(i_0)=t_1$ and $s(j_0)=t_r$. 
\item[(c)] All other vertices and edges come in pairs of two $\sigma$-conjugates.
\item[(d)] Let $T/\gensigma$ denote the graph whose set of vertices is $V(T)/\gensigma$, and where two vertices $\{t,\sigma(t)\}$ and $\{t',\sigma(t')\}$ are  joined by an edge if and only if $t$ is joined by an edge to $t'$ or to $\sigma(t')$.
Then the dual tree of $\barX_0 = X_0/\gensigma$ is naturally isomorphic to $T/\gensigma$.
\item[(e)] Any subtree of $T$ which contains $t_1,\ldots,t_r$ and at most one additional edge emanating from each of these vertices maps isomorphically to its image in $T/\gensigma$.
\end{enumerate}
\end{Prop}

\begin{Proof}
Assertions (a) through (c) follow from the construction and Proposition \ref{SigmaActProp}. Assertion (d) follows from the definition of $X_0/\gensigma$ and Proposition \ref{QuotientProp}. Note that for an arbitrary finite tree with an automorphism of order $2$ the quotient described in (d) is not necessarily a tree, but in this case that follows from assertions (a) to (c). Finally, assertion (e) is a direct consequence of (a) and (b).
\end{Proof}

\medskip
The sketch below shows what $T$ and $T/\gensigma$ might typically look like:
\begin{center}
\parbox[t]{.4\textwidth}
{\setlength{\fboxsep}{10pt}
\fbox{\xy 
% (-30,5)*{}; 
% (30,60)*{}; 
(0,64)*{\bullet};
(0,56)*{\bullet};
(0,40)*{\bullet};
(0,32)*{\bullet};
(0,63)*{};(0,57)*{};**\crv{};
(0,39)*{};(0,33)*{};**\crv{};
(-7,56)*{};(7,56)*{};**\crv{};
(.7,55)*{};(.7,41)*{};**\crv{~*=<4pt>{}~**=<4pt>{.}};
(8,56)*{\bullet};
(8,36)*{\bullet};
(8,28)*{\bullet};
(8,20)*{\bullet};
(8,12)*{\bullet};
(12,4)*{\bullet};
(16,8)*{\bullet};
(12,48)*{\bullet};
(20,44)*{\bullet};
(28,48)*{\bullet};
(28,44)*{\bullet};
(28,40)*{\bullet};
(8,56)*{};(12,48)*{};**\crv{};
(12,48)*{};(20,44)*{};**\crv{};
(20,44)*{};(28,48)*{};**\crv{};
(21,44)*{};(27,44)*{};**\crv{};
(20,44)*{};(28,40)*{};**\crv{};
(0,32)*{};(8,36)*{};**\crv{};
(0,32)*{};(8,28)*{};**\crv{};
(8,27)*{};(8,13)*{};**\crv{};
(8,12)*{};(12,4)*{};**\crv{};
(8,12)*{};(16,8)*{};**\crv{};
(-8,56)*{\bullet};
(-8,36)*{\bullet};
(-8,28)*{\bullet};
(-8,20)*{\bullet};
(-8,12)*{\bullet};
(-12,4)*{\bullet};
(-16,8)*{\bullet};
(-12,48)*{\bullet};
(-20,44)*{\bullet};
(-28,48)*{\bullet};
(-28,44)*{\bullet};
(-28,40)*{\bullet};
(-8,56)*{};(-12,48)*{};**\crv{};
(-12,48)*{};(-20,44)*{};**\crv{};
(-20,44)*{};(-28,48)*{};**\crv{};
(-21,44)*{};(-27,44)*{};**\crv{};
(-20,44)*{};(-28,40)*{};**\crv{};
(0,32)*{};(-8,36)*{};**\crv{};
(0,32)*{};(-8,28)*{};**\crv{};
(-8,27)*{};(-8,13)*{};**\crv{};
(-8,12)*{};(-12,4)*{};**\crv{};
(-8,12)*{};(-16,8)*{};**\crv{};
(3,64)*{\scriptstyle t_1};
(3,54)*{\scriptstyle t_2};
(4.5,40)*{\scriptstyle t_{r-1}};
(2,29)*{\scriptstyle t_r};
\endxy}
\\[6pt]
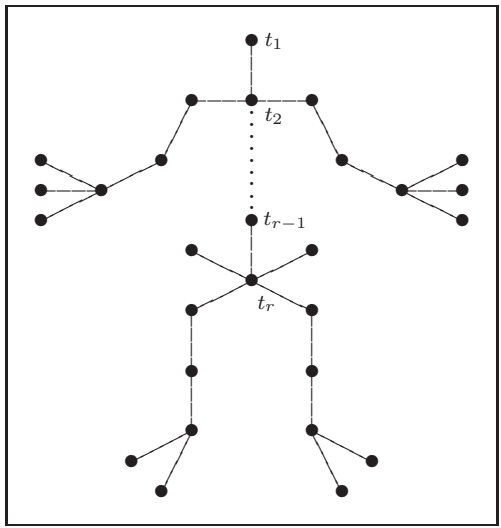
Dual tree $T$, where $\sigma$ is the reflection
at the vertical axis of symmetry}
\qquad\qquad  
\parbox[t]{.27\textwidth}
{\setlength{\fboxsep}{10pt}
\fbox{\xy 
% (-30,5)*{}; 
% (30,60)*{}; 
(0,64)*{\bullet};
(0,56)*{\bullet};
(0,40)*{\bullet};
(0,32)*{\bullet};
(0,63)*{};(0,57)*{};**\crv{};
(0,39)*{};(0,33)*{};**\crv{};
(1,56)*{};(7,56)*{};**\crv{};
(.7,55)*{};(.7,41)*{};**\crv{~*=<4pt>{}~**=<4pt>{.}};
(8,56)*{\bullet};
(8,36)*{\bullet};
(8,28)*{\bullet};
(8,20)*{\bullet};
(8,12)*{\bullet};
(12,4)*{\bullet};
(16,8)*{\bullet};
(12,48)*{\bullet};
(20,44)*{\bullet};
(28,48)*{\bullet};
(28,44)*{\bullet};
(28,40)*{\bullet};
(8,56)*{};(12,48)*{};**\crv{};
(12,48)*{};(20,44)*{};**\crv{};
(20,44)*{};(28,48)*{};**\crv{};
(21,44)*{};(27,44)*{};**\crv{};
(20,44)*{};(28,40)*{};**\crv{};
(0,32)*{};(8,36)*{};**\crv{};
(0,32)*{};(8,28)*{};**\crv{};
(8,27)*{};(8,13)*{};**\crv{};
(8,12)*{};(12,4)*{};**\crv{};
(8,12)*{};(16,8)*{};**\crv{};
(-3,64)*{\scriptstyle \bar t_1};
(-3,56)*{\scriptstyle \bar t_2};
(-4.5,40)*{\scriptstyle \bar t_{r-1}};
(-3,32)*{\scriptstyle \bar t_r};
\endxy}
\\[6pt]
The quotient $T/\gensigma$}
\end{center}

%%%%%%%%%%%%%%%%%%%%%%%%%%%%%%%%%%%%%%%%%%%%%%%%%%%%%%%%%%%%%%%%%%%%%%%%%%%%%%%%%%%%%%%%%%%%%

\section{Stable curves associated to quadratic morphisms}
\label{StableCurvesQuad}

Now we fix a discrete valuation ring $R$ with quotient field $K$ and residue field $k$ of characteristic $\not=2$. We impose no other restrictions on the characteristics of $K$ and~$k$. We fix a finite mapping scheme $(\Gamma,\tau,i_1,j_1)$ with $|\Gamma|\ge3$. We also fix a $\Gamma$-marked quadratic morphism $(C,f,P,Q,s)$ over~$K$. Without loss of generality we may assume that $C=\BP^1_K$ and $P=0$ and $Q=\infty$. Then $f$ has the form $f(x) = \frac{ax^2+b}{cx^2+d}$ for some $\binom{\,a\ \ b\,}{c\ \ d}\in\PGL_2(K)$ and is a Galois covering of degree $2$ with the non-trivial covering automorphism $\sigma\colon x\mapsto -x$. The critical points $0$ and $\infty$ are then precisely the fixed points of~$\sigma$.

\medskip
To analyze the possible degeneration of $(C,f,P,Q,s)$ over $\Spec R$ we view $C$ as a smooth marked curve of genus zero in two ways. The first is simply $C$ with the marking $s\colon \Gamma\into C(K)$.

\medskip
For the second recall from Definition \ref{AbsPostCritDef} (c) that there exists at most one element $i_0\in\Gamma$ with $\tau(i_0)=i_1$. If there exists none, we choose a new symbol $i_0\not\in\Gamma$. Likewise, there exists at most one element $j_0\in\Gamma$ with $\tau(j_0)=j_1$, and if there is none, we choose another new symbol $j_0\not\in\Gamma\cup\{i_0\}$. In either case we set $\Gamma_0 := \Gamma\cup\{i_0,j_0\}$ and define
\UseTheoremCounterForNextEquation
\begin{equation}\label{Gamma0Cons}
\tilde\Gamma\ :=\ \Gamma_0 \;\sqcup\; 
\bigl\{\sigma(\gamma)\,\bigm|\,\gamma\in\Gamma, \ 
\gamma\not=i_0,j_0,\ 
|\tau^{-1}(\tau(\gamma))| =1 \bigr\},
\end{equation}
where the $\sigma(\gamma)$ are distinct new symbols $\not\in\Gamma_0$.
We define a map $\sigma\colon\tilde\Gamma\to\tilde\Gamma$ by setting
\UseTheoremCounterForNextEquation
\begin{equation}\label{SigmaCons}
\left\{\begin{array}{ll}
\gamma\mapsto\gamma & \hbox{for $\gamma\in\{i_0,j_0\}$,}\\[3pt]
\gamma\mapsto\gamma'\mapsto\gamma & \hbox{for all $\gamma\in\Gamma$ with $\tau^{-1}(\tau(\gamma)) = \{\gamma,\gamma'\}$ and $\gamma\not=\gamma'$,} \\[3pt]
\gamma\mapsto\sigma(\gamma)\mapsto\gamma & \hbox{for all other $\gamma\in\Gamma$.}
\end{array}\right.
\end{equation}
By construction this is an automorphism of~$\tilde\Gamma$ of order $2$ with precisely two fixed points $i_0$ and~$j_0$. We also define a map $\tilde\tau\colon \tilde\Gamma\to\Gamma$ by setting
\UseTheoremCounterForNextEquation
\begin{equation}\label{TauTildeCons}
\left\{\begin{array}{ll}
i_0\mapsto i_1 & \\[3pt]
j_0\mapsto j_1 & \\[3pt]
\gamma\mapsto \tau(\gamma) & \hbox{for all $\gamma\in\Gamma$,}\\[3pt]
\sigma(\gamma)\mapsto \tau(\gamma) & \hbox{for all $\sigma(\gamma)\in\tilde\Gamma\setminus\Gamma_0$,}
%\gamma\mapsto\sigma(\gamma)\mapsto\gamma & \hbox{for all other $\gamma\in\Gamma$.}
\end{array}\right.
\end{equation}
which is well-defined by the definition of $i_0$ and~$j_0$. Definition \ref{AbsPostCritDef} (a) shows that $\tilde\tau$ is surjective, and Definition \ref{AbsPostCritDef} (b) implies that for any $\gamma$, $\gamma'\in\tilde\Gamma$ we have $\tilde\tau(\gamma) = \tilde\tau(\gamma')$ if and only if $\gamma'=\gamma$ or $\sigma(\gamma)$. Thus $\tilde\tau$ induces a bijection
\UseTheoremCounterForNextEquation
\begin{equation}\label{GammaQuot}
\tilde\Gamma/\gensigma \stackrel{\sim}{\longto} \Gamma.
\end{equation}
The composite $\Gamma\into\tilde\Gamma\onto\tilde\Gamma/\gensigma \stackrel{\sim}{\to} \Gamma$ is, of course, just the map~$\tau$.
Furthermore, observe that if $i_0\in\Gamma$, then Definition \ref{GammaMarkDef} implies that $f(s(i_0)) = s(\tau(i_0)) = s(i_1) = f(P)$, which by (\ref{RamLift}) implies that $s(i_0)=P$. The same argument shows that, if $j_0\in\Gamma$, then $s(j_0)=Q$. Thus the following map $\tilde s\colon \tilde\Gamma\to C(K)$ is well-defined:
\UseTheoremCounterForNextEquation
\begin{equation}\label{SigmaTildeCons}
\left\{\begin{array}{ll}
i_0 \mapsto P & \\[3pt]
j_0 \mapsto Q & \\[3pt]
\gamma\mapsto s(\gamma) & \hbox{for all $\gamma\in\Gamma$,}\\[3pt]
\sigma(\gamma) \mapsto \sigma(s(\gamma)) & \hbox{for all $\sigma(\gamma)\in\tilde\Gamma\setminus\Gamma_0$.}
\end{array}\right.
\end{equation}
{}From the construction of $\tilde\Gamma$ and $\tilde s$ we readily deduce that $\tilde s$ is injective. In fact, the image $\tilde s(\Gamma_0)$ is precisely the 
non-strictly-postcritical orbit ${\{ f^n(P), f^n(Q) \mid n\ge0 \}}$ of~$f$. Together the construction implies that the following diagram commutes:
\UseTheoremCounterForNextEquation
\begin{equation}\label{SPlusSTilde}
\vcenter{\xymatrix{
C(K) \ar@{=}[r] & C(K) \ar[r]^\sigma & C(K) \ar[r]^f & C(K) \\
\Gamma \ar@{^{ (}->}[r] \ar@{^{ (}->}[u]^{s} 
& \tilde\Gamma \ar[r]^\sigma \ar@{^{ (}->}[u]^{\tilde s} 
& \tilde\Gamma \ar@{->>}[r]^{\tilde\tau} \ar@{^{ (}->}[u]^{\tilde s} 
& \Gamma\rlap{.} \ar@{^{ (}->}[u]^{s} \\}}
\end{equation}

\medskip
Now recall that $|\Gamma|\ge3$ and hence $|\tilde\Gamma|\ge3$. Thus $(C,s)$ is a smooth stable $\Gamma$-marked curve of genus zero, and $(C,\tilde s)$ is a smooth stable $\tilde\Gamma$-marked curve of genus zero.
Let $(X,s)$ denote the stable marked curve of genus zero over $\Spec R$ with disjoint sections $s(\gamma)$ for all $\gamma\in\Gamma$ which extends $(C,s)$.
Let $(\tilde X,\tilde s)$ denote the stable marked curve of genus zero over $\Spec R$ with disjoint sections $\tilde s(\tilde\gamma)$ for all $\tilde\gamma\in\tilde\Gamma$ which extends $(C,\tilde s)$. 
Since $s=\tilde s|\Gamma$ over~$K$, the identity on $C$ extends to a morphism $\tilde X\onto X$, and the pair $(X,s)$ can be obtained from $(\tilde X,\tilde s|\Gamma)$ by contracting suitable irreducible components of the special fiber, as explained in Section~\ref{StableCurves}.
Also, by the uniqueness of $(\tilde X,\tilde s)$, the involution $\sigma$ of $C$ extends to an automorphism of $\tilde X$ which still satisfies $\sigma^2=\id$ and $\sigma\circ\tilde s = \tilde s\circ\sigma$. 

\begin{Prop}\label{FExtends}
The morphism $f$ extends to a unique morphism $f\colon\tilde X\to X$ and induces an isomorphism 
$$\tilde X/\gensigma \stackrel{\sim}{\longto} X.$$
\end{Prop}

\begin{Proof}
As in Section \ref{StableAut} the quotient $\tilde X/\gensigma$ inherits a marking indexed by $\tilde\Gamma/\gensigma$ and becomes a stable marked curve over $\Spec R$ by Proposition \ref{QuotientProp} (a). In the generic fiber of $X$ the involution $\sigma$ is precisely the non-trivial covering automorphism of~$f$, so that $f$ induces an isomorphism between the generic fibers of $\tilde X/\gensigma$ and~$X$. This isomorphism is compatible with the markings via the isomorphism $\tilde\Gamma/\gensigma \stackrel{\sim}{\to} \Gamma$ from (\ref{GammaQuot}). Since the stable extension of a marked curve is unique up to unique isomorphism, it follows that the isomorphism extends to an isomorphism $\tilde X/\gensigma \stackrel{\sim}{\to} X$.
\end{Proof}

\medskip
{}From the commutative diagram (\ref{SPlusSTilde}) and Proposition \ref{FExtends} we deduce the commutative diagram
\UseTheoremCounterForNextEquation
\begin{equation}\label{SPlusSTildeX}
\vcenter{\xymatrix{
X(R) & \tilde X(R) \ar[l] \ar[r]^\sigma & \tilde X(R) \ar[r]^{\smash{f}} & X(R) \\
\Gamma \ar@{^{ (}->}[r] \ar@{^{ (}->}[u]^{s} 
& \tilde\Gamma \ar[r]^\sigma \ar@{^{ (}->}[u]^{\tilde s} 
& \tilde\Gamma \ar@{->>}[r]^{\tilde\tau} \ar@{^{ (}->}[u]^{\tilde s} 
& \Gamma\rlap{.} \ar@{^{ (}->}[u]^{s} \\}}
\end{equation}

\medskip
Now we translate this information into combinatorial information on dual trees. 
% Let $X_0$ and $\tilde X_0$ denote the respective closed fibers of $X$ and~$\tilde X$, and let $T$ and $\tilde T$ denote their respective dual trees. 
Let $T$ and $\tilde T$ denote the respective dual trees of the closed fibers of $X$ and~$\tilde X$. As in (\ref{MarkOnTree}) we obtain natural maps $s\colon \Gamma\to V(T)$ and $\tilde s\colon \tilde\Gamma\to V(\tilde T)$. Combining (\ref{StabilizationDiag}) with Propositions \ref{SigmaTreeProp} (d) and \ref{FExtends} we obtain a diagram 
\UseTheoremCounterForNextEquation
\begin{equation}\label{SPlusSTildeT}
\vcenter{\xymatrix{
V(T) \ar@{^{ (}->}[r] \ar@{}[dr]|{\textstyle!!!}
& V(\tilde T) \ar[r]^\sigma & V(\tilde T) \ar@{->>}[r] & 
V(\tilde T)/\gensigma \ar[r]^\sim & V(T) \\
\Gamma \ar@{^{ (}->}[r] \ar@{^{ (}->}[u]^{s} 
& \tilde\Gamma \ar[r]^\sigma \ar@{^{ (}->}[u]^{\tilde s} 
& \tilde\Gamma \ar@{->>}[r] \ar@{^{ (}->}[u]^{\tilde s} 
\ar@/_15pt/[rr]_{\tilde\tau}
& \tilde\Gamma/\gensigma \ar@{^{ (}->}[u] \ar[r]^{\sim}
& \Gamma\rlap{.} \ar@{^{ (}->}[u]^{s} \\}}
\end{equation}
This commutes everywhere except at the leftmost square marked !!!, for which the modified rule is explained in Proposition \ref{StabilizationProp} (c). The ingredients for the next section are the combinatorics of this diagram and the local description of the quotient from Proposition \ref{QuotientProp}~(b).

We will abbreviate the marked vertices of $T$ and $\tilde T$ as follows:
$$\left\{\!\begin{array}{c}
P_n := s(i_n) \\[3pt]
Q_n := s(j_n))
\end{array}\!\right\}
\quad\hbox{for all $n\ge1$;}\qquad
\left\{\!\begin{array}{c}
\tilde P_n := \tilde s(i_n)) \\[3pt]
\tilde Q_n := \tilde s(j_n))
\end{array}\!\right\}
\quad\hbox{for all $n\ge0$.}$$
Under the inclusion $V(T)\into V(\tilde T)$ at the left hand side of the diagram (\ref{SPlusSTildeT}), the vertices $P_n$ and $Q_n$ can be constructed from $\tilde P_n$ and $\tilde Q_n$ by the procedure in Proposition \ref{StabilizationProp} (c).
The map $V(\tilde T)\onto V(T)$ on the right half of the diagram (\ref{SPlusSTildeT}) sends $\tilde P_n$ to $P_{n+1}$ and $\tilde Q_n$ to~$Q_{n+1}$.
Also recall from Proposition \ref{SigmaTreeProp} (a) that the fixed points of $\sigma$ on $V(\tilde T)$ are precisely the vertices between and including $\tilde P_0 = t_1 = \tilde s(i_0)$ and $\tilde Q_0 = t_r = \tilde s(j_0)$.

\begin{Lem}\label{SpineSurvives}
\begin{enumerate}
\item[(a)] Any vertex strictly between $\tilde P_0$ and $\tilde Q_0$ survives in~$T$.
\item[(b)] If $\tilde P_0\not=\tilde Q_0$, then $\tilde P_0$ survives in $T$ unless one of the following happens in~$\tilde T$:
\begin{enumerate}
\item[(i)] $i_0\not\in\Gamma$, there is only one edge emanating from~$\tilde P_0$, and the only markings at $\tilde P_0$ are $\tilde P_0=s(i_0)=s(\gamma)=\sigma(s(\gamma))$ for a unique $\gamma\in\Gamma$ with $\sigma(\gamma)\not\in\Gamma$. 
\item[(ii)] $i_0\not\in\Gamma$, there is no other marking at~$\tilde P_0$, and the connected components of $\tilde T\setminus\{\tilde P_0\}$ are precisely that containing~$\tilde Q_0$ and two others of the form $S$ and~$\sigma(S)$, where $\sigma(S)$ possesses no markings of the form $s(\gamma)$ for $\gamma\in\Gamma$.
\end{enumerate}
\end{enumerate}
\end{Lem}

\begin{Proof}
Let $t$ be any vertex between or equal to $\tilde P_0$ and~$\tilde Q_0$. By Proposition \ref{StabilizationProp} (a) it survives in $T$ if and only if 
\UseTheoremCounterForNextEquation
\begin{equation}\label{Survfor}
\bigl|\tilde s^{-1}(t)\cap \Gamma\bigr| + \left|\left\{
\begin{array}{c}
\hbox{$S$ connected component of $\tilde T\setminus\{t\}$} \\
\hbox{such that $\tilde s^{-1}(V(S))\cap\Gamma\not=\emptyset$}
\end{array}
\right\}\right| \ge3.
\end{equation}
First note that, by the stability of~$\tilde T$, the vertex $t$ must either have another marking or another edge emanating from it, or both. In the second case the connected component of $\tilde T\setminus\{t\}$ determined by that edge must again contain some mark by stability. Moreover, the marking at~$t$ or in that connected component must correspond to an element of $\tilde\Gamma$ that is not fixed by~$\sigma$, and hence have the form $s(\gamma)$ or $\sigma(s(\gamma))$ for some $\gamma\in\Gamma$. Since $t$ itself is fixed by~$\sigma$, we may conjugate the connected component by~$\sigma$, if necessary, after which the marking is $s(\gamma)$ for some $\gamma\in\Gamma$. Then in either case, this marked point contributes $1$ to the left hand side of (\ref{Survfor}).

If $t$ lies strictly between $\tilde P_0$ and~$\tilde Q_0$, we can apply the preceding remarks also to $\tilde P_0$ and~$\tilde Q_0$ in place of~$t$. It follows that the two connected components of $\tilde T\setminus\{t\}$ containing $\tilde P_0$, respectively~$\tilde Q_0$, each contribute $1$ to the left hand side of (\ref{Survfor}). Thus the total sum is $\ge3$, and so $t$ survives in~$T$, as desired.

Now assume that $t=\tilde P_0\not=\tilde Q_0$. Then the connected component of $\tilde T\setminus\{t\}$ containing $\tilde Q_0$ contributes $1$ to the left hand side of (\ref{Survfor}). By the above remarks there must be another marking at $t$ or another connected component of $\tilde T\setminus\{t\}$ that also contributes~$1$. Thus the only way in which $t$ cannot survive is that there is no other contribution. This requires first of all that the marking $\tilde P_0 = \tilde s(i_0)$ itself is removed, i.e., that $i_0\not\in\Gamma$. Next it requires that, if there is another marking at~$t$, then $t=s(\gamma)=\sigma(s(\gamma))$ for a unique 
$\gamma\in\Gamma$ with $\sigma(\gamma)\not\in\Gamma$, and there is no other edge emanating from~$t$. This is the case (i) of the lemma. Thirdly, if there is no other marking at~$t$, there must be a unique connected component $S$ of $\tilde T\setminus\{t\}$ which contributes $1$ to the left hand side of (\ref{Survfor}). Then $\tilde T\setminus\{t\}$ consists of the three connected components $S$ and $\sigma(S)$ and that containing~$\tilde Q_0$, where $\sigma(S)$ cannot contain any marking of the form $s(\gamma)$ for $\gamma\in\Gamma$. This is the case (ii) of the lemma, and we are done.
\end{Proof}

%%%%%%%%%%%%%%%%%%%%%%%%%%%%%%%%%%%%%%%%%%%%%%

\begin{Ex}\label{NotAllOccur2}
\rm We continue the example in Remark \ref{NotAllOccur} with the mapping scheme $\Gamma$
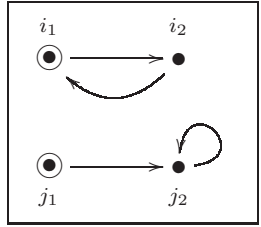
$$\fbox{\ $\xymatrix@R-25pt@C-10pt{
\scriptstyle i_1 && \scriptstyle i_2 &\\ 
\circbullet\ar[rr]&&\bullet\ar@/^15pt/[ll] & \\
{\phantom{X}}&&&\\
{\phantom{X}}&&&\\
\circbullet\ar[rr]&&\bullet \ar@(r,u)[] & \\
\vphantom{x}\mysmash{2pt}{j_1} && \mysmash{2pt}{j_2} & \\ 
}$}\qquad$$
The construction of $\tilde\Gamma$ implies that in this case $i_0=i_2$ and
$$\tilde \Gamma = \bigl\{i_1,i_2,\sigma(i_1),j_0,j_1,j_2\bigr\}.$$
Consider the quadratic morphism $(\BP^1_\BQ, {x\mapsto \frac{x^2-4}{x^2+2}},0,\infty)$ over $K:=\BQ$, which has postcritical orbit~$\Gamma$. The induced $\tilde\Gamma$-marking $\tilde s$ on $\BP^1_\BQ$ is given by 
$$\left\{\!\begin{array}{rcc}
i_0=i_2 &\!\!\mapsto\!\!\!& 0 \\[3pt]
i_1 &\!\!\mapsto\!\!\!& -2 \\[3pt]
\sigma(i_1) &\!\!\mapsto\!\!\!& 2 \\[3pt]
\end{array}\!\right\}
\qquad\hbox{and}\qquad
\left\{\!\begin{array}{rcc}
j_0 &\!\!\mapsto\!\!\!& \infty \\[3pt]
j_1 &\!\!\mapsto\!\!\!& 1 \\[3pt]
j_2 &\!\!\mapsto\!\!\!& -1
\end{array}\!\right\}_.$$
Here the images of $i_1$ and $j_1$, respectively of $\sigma(i_1)$ and $j_2$, are congruent modulo~$3$, but there are no other congruences modulo~$3$. Thus the stable extension $\tilde X$ over $R:=\BZ_{(3)}$ is obtained by blowing up $\BP^1\times\Spec\BZ_{(3)}$ in the two points $x=\pm1$ of the special fiber. The special fiber of $\tilde X$ has three irreducible components, one met by the sections $0$ and~$\infty$, one by the sections $-2$ and~$1$, and one by the sections $2$ and~$-1$, and the first of these irreducible components meets each of the other two in an ordinary double point.

The stable extension $X$ can be obtained from $\tilde X$ by removing the sections $\tilde s(j_0)$ and $\tilde s(\sigma(i_1))$ and contracting unstable irreducible components of the special fiber. Here the blowup at $x=-1$ in the special fiber is reversed, while the blowup at $x=1$ in the special fiber remains, because the points $i_1$ and $j_1$ remain in~$\Gamma$. The closed fiber of $X$ possesses two irreducible components, one met by the sections $0$ and~$-1$, the other by the sections $-2$ and~$1$. One easily checks that the same description of the closed fiber of $X$ is obtained from the isomorphism $\tilde X/\gensigma\cong X$ of Proposition \ref{FExtends}.

The associated dual trees with their markings and the maps between them are sketched below:
$$
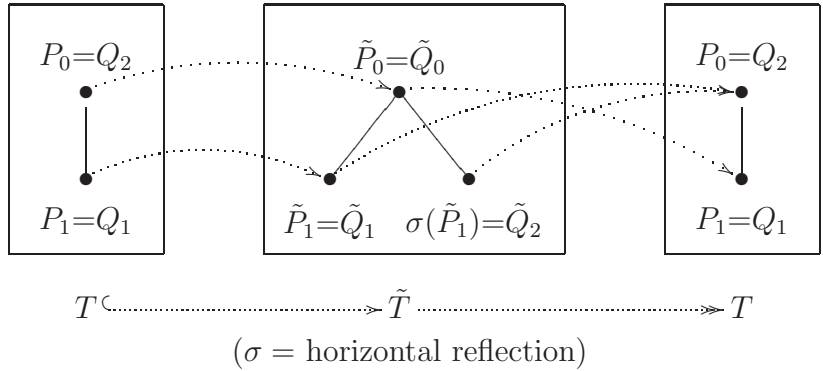
\xymatrix@C=5pt@R=5pt{
&&&&&&&&&&&&&&&&&&&&&&&&& \\
&&&&&&&&&&&&&&&&&&&&&&&&& \\
&&&&&&&&&&&&&&&&&&&&&&&&& \\
&& \mybullet \ar@{-}[ddd] \ar@{}[uu]|>>{\textstyle P_0{=}Q_2}
\ar@{.>}@/^10pt/@<1pt>[rrrrrrrrrr(0.97)] &&&&&&&&&& 
\mybullet \ar@{}[uu]|>>{\textstyle\tilde P_0{=}\tilde Q_0}
% \ar@{-}[dddll]\ar@{-}[dddrr] 
\ar@{.>}@/^10pt/@<1pt>[rrrrrrrrrrrddd(0.97)] &&&&&&&&&&& 
\mybullet \ar@{-}[ddd] \ar@{}[uu]|>>{\textstyle P_0{=}Q_2} && \\
&&&&&&&&&&&&&&&&&&&&&&&&& \\
&&&&&&&&&&&&&&&&&&&&&&&&& \\
&& \mybullet \ar@{}[dd]|>>>{\textstyle P_1{=}Q_1}
\ar@{.>}@/^10pt/@<1pt>[rrrrrrrr(0.96)] &&&&&&&& 
\mybullet \ar@{.>}@/^15pt/@<.8pt>[rrrrrrrrrrrrruuu(0.976)] 
\ar@{-}[uuurr(1.03)] 
\ar@{}[dd]|>>>{\textstyle\tilde P_1{=}\tilde Q_1} &&&& 
\mybullet \ar@{.>}@/^10pt/@<1pt>[rrrrrrrrruuu(0.97)] 
\ar@{-}[uuull(1.03)] 
\ar@{}[dd]|>>>{\textstyle\ \sigma(\tilde P_1){=}\tilde Q_2} &&&&&&&&& 
\mybullet  \ar@{}[dd]|>>>{\textstyle P_1{=}Q_1} && \\
&&&&&&&&&&&&&&&&&&&&&&&&& \\
&&&&&&&&&&&&&&&&&&&&&&&&& \\
&&&&&&&&&&&&&&&&&&&&&&&&& \\
\save "1,1"."9,5"*[F]\frm{}\restore
\save "1,9"."9,18"*[F]\frm{}\restore
\save "1,22"."9,26"*[F]\frm{}\restore
&& T \ar@{^{ (}.>}[rrrrrrrrrr] &&&&&&&&&& \tilde T \ar@{.>>}[rrrrrrrrrrr] &&&&&&&&&&& T && \\
&&&&&&&&&&&&
\save[]*\txt{\ \ ($\sigma=$ horizontal reflection)}\restore
&&&&&&&&&&& 
&& \\
}$$
\end{Ex}

\begin{Ques}\label{WhichOccur}
\rm As some non-trivial types of degeneration can occur, while in the next section we will exclude others, it seems interesting to ask which types of degeneration are actually possible in general.
\end{Ques}

%%%%%%%%%%%%%%%%%%%%%%%%%%%%%%%%%%%%%%%%%%%%%%%%%%%%%%%%%%%%%%%%%%%%%%%%%%%%%%%%%%%%%%%%%%%%%

\section{Excluding certain types of bad reduction}
\label{Exclude}

% Recall that a subtree of a tree is the graph obtained by deleting some vertices and all edges adjacent to deleted vertices, provided that the result is still a tree. 
%In this section we prove that $\tilde T$ does not contain subtrees of certain kinds with certain marked points. In each case we indicate a subtree graphically by letting a solid line denote an edge, and a dashed line a connected sequence of distinct edges, possibly empty. Thus two vertices connected by a solid line are distinct, but vertices connected by one or more dashed lines may be equal. Sometimes we label a dashed line by the length of the associated path, i.e., by the number of its edges. Labeling a dashed line by `$\,>0\,$' means that its length is $>0$, and in that case vertices separated by it are distinct.
Keeping the notation of the preceding section, we now prove that $\tilde T$ does not contain subtrees of certain kinds with certain marked points. In each case we indicate a subtree graphically by letting a dashed line denote a connected sequence of distinct edges, possibly empty. Thus vertices connected by one or more dashed lines may be equal. Sometimes we label a dashed line by the length of the associated path, i.e., by the number of its edges. Labeling a dashed line by `$\,>0\,$' means that its length is $>0$, and in that case vertices separated by it are distinct.

We begin with the following lemma, whose proof is typical for the arguments in this section.

\begin{Lem}\label{An}
The tree $\tilde T$ does not contain a subtree of the form:
$$\fbox
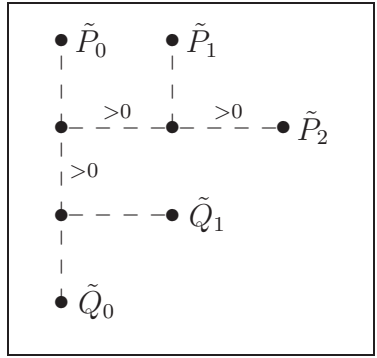
{\xymatrix@C=8pt@R=5pt{
&&&&&&&& \\
& \mybullet \ar@{--}[ddd] \ar@{}[r]|-{\ \ \textstyle\tilde P_0} &&& 
\mybullet \ar@{--}[ddd] \ar@{}[r]|-{\ \ \textstyle\tilde P_1} &&&& \\
&&&&&&&& \\
&&&&&&&& \\
& \mybullet \ar@{--}[ddd]^{>0} \ar@{--}[rrr]^{>0} &&& 
\mybullet \ar@{--}[rrr]^{>0} &&&
\mybullet \ar@{}[r]|-{\ \ \textstyle\tilde P_2} & \\
&&&&&&&& \\
&&&&&&&& \\
& \mybullet \ar@{--}[ddd]  \ar@{--}[rrr] &&& 
\mybullet \ar@{}[r]|-{\ \ \ \textstyle\tilde Q_1} &&&& \\
&&&&&&&& \\
&&&&&&&& \\
& \mybullet \ar@{}[r]|-{\ \ \ \textstyle\tilde Q_0} &&&&&&& \\
&&&&&&&&& \\}}$$
\end{Lem}

\begin{Proof}
For any $n\ge2$ let $A_n$ denote a subtree of $\tilde T$ of the form:
\UseTheoremCounterForNextEquation
\begin{equation}\label{An1}
\fbox{\xymatrix@C=8pt@R=5pt{
&&&&&&&&&&&&&&&&& \\
& \mybullet \ar@{--}[ddd] \ar@{}[r]|-{\ \ \textstyle\tilde P_0} &&& 
\mybullet \ar@{--}[ddd] \ar@{}[r]|-{\ \ \textstyle\tilde P_1} &&&&&&&&
\mybullet \ar@{--}[ddd] \ar@{}[r]|-{\ \ \ \ \;\textstyle\tilde P_{n-1}} &&&&& \\
&&&&&&&&&&&&&&&&& \\
&&&&&&&&&&&&&&&&& \\
& \mybullet \ar@{--}[ddd]^{>0} \ar@{--}[rrr]^{>0} &&& 
\mybullet \ar@{--}[rrr]^{>0} \ar@{}[ddr]|<<<<{\textstyle t} &&& 
\ar@{.}[rr] && \ar@{--}[rrr]^{>0} &&& 
\mybullet \ar@{--}[rrr]^{>0} &&& 
\mybullet \ar@{}[r]|-{\ \ \textstyle\tilde P_n} && \\
&&&&&&&&&&&&&&&&& \\
&&&&&&&&&&&&&&&&& \\
& \mybullet \ar@{--}[ddd]  \ar@{--}[rrr] &&& 
\mybullet \ar@{}[r]|-{\ \ \ \textstyle\tilde Q_1} &&&&&&&&&&&&& \\
&&&&&&&&&&&&&&&&& \\
&&&&&&&&&&&&&&&&& \\
& \mybullet \ar@{}[r]|-{\ \ \ \textstyle\tilde Q_0} &&&&&&&&&&&&&&&& \\
&&&&&&&&&&&&&&&&& \\}}
\end{equation}
We must show that $\tilde T$ does not contain a subtree of the form~$A_2$. 
Suppose to the contrary that it does. We will prove by induction that $\tilde T$ then contains a subtree of the form $A_n$ for every $n\ge2$. But in a subtree of the form $A_n$ the vertices $\tilde P_0,\ldots,\tilde P_n$ are all distinct, which is impossible for $n\ge|\Gamma|$. Thus we obtain a contradiction, and the lemma follows.

For the induction step assume that $\tilde T$ contains a subtree of the form $A_n$ for some $n\ge2$. In this subtree, any vertex that is fixed by~$\sigma$, that is, any vertex equal to $\tilde P_0$ or $\tilde Q_0$ or between them, has at most one neighbor that is not fixed by~$\sigma$. By Proposition \ref{SigmaTreeProp} (e) this implies that the subtree maps isomorphically to its image in $\tilde T/\gensigma$. In view of the right half of the diagram (\ref{SPlusSTildeT}) it follows that $T$ contains a subtree of the form:
\UseTheoremCounterForNextEquation
\begin{equation}\label{An1t}
\fbox{\xymatrix@C=8pt@R=5pt{
&&&&&&&&&&&&&&&&& \\
& \mybullet \ar@{--}[ddd] \ar@{}[r]|-{\ \ \ \textstyle P_1} &&& 
\mybullet \ar@{--}[ddd] \ar@{}[r]|-{\ \ \ \textstyle P_2} &&&&&&&&
\mybullet \ar@{--}[ddd] \ar@{}[r]|-{\ \ \ \textstyle P_n} &&&&& \\
&&&&&&&&&&&&&&&&& \\
&&&&&&&&&&&&&&&&& \\
& \mybullet \ar@{--}[ddd]^{>0} \ar@{--}[rrr]^{>0} &&& 
\mybullet \ar@{--}[rrr]^{>0} &&&
\ar@{.}[rr] && \ar@{--}[rrr]^{>0} &&& 
\mybullet \ar@{--}[rrr]^{>0} &&& 
\mybullet \ar@{}[r]|-{\ \ \ \ \;\textstyle P_{n+1}} && \\
&&&&&&&&&&&&&&&&& \\
&&&&&&&&&&&&&&&&& \\
& \mybullet \ar@{}[r]|-{\ \ \ \ \textstyle Q_1} &&&&&&&&&&&&&&&& \\
&&&&&&&&&&&&&&&&& \\}}
\end{equation}
On the other hand $T$ can be obtained from $\tilde T$ by the stabilization procedure described in Section \ref{StableCurves}. Conversely this means that $\tilde T$ can be obtained from $T$ by inserting edges and thereby moving certain marked points apart or further apart, but never the opposite. Thus from (\ref{An1t}) it follows that $\tilde T$ contains a subtree of the form:
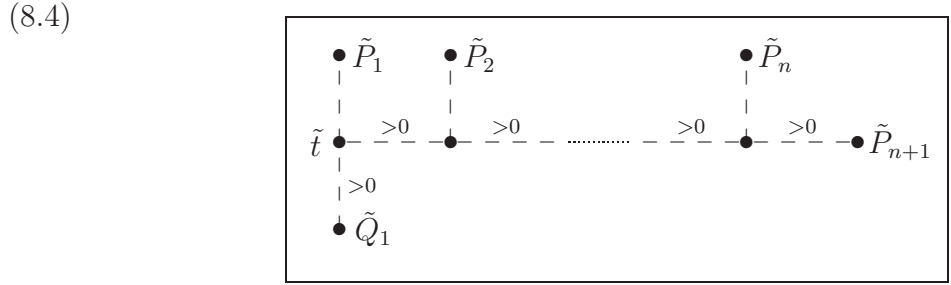
\UseTheoremCounterForNextEquation
\begin{equation}\label{An2}
\fbox{\xymatrix@C=8pt@R=5pt{
&&&&&&&&&&&&&&&&& \\
& \mybullet \ar@{--}[ddd] \ar@{}[r]|-{\ \ \textstyle\tilde P_1} &&& 
\mybullet \ar@{--}[ddd] \ar@{}[r]|-{\ \ \textstyle\tilde P_2} &&&&&&&&
\mybullet \ar@{--}[ddd] \ar@{}[r]|-{\ \ \textstyle\tilde P_n} &&&&& \\
&&&&&&&&&&&&&&&&& \\
&&&&&&&&&&&&&&&&& \\
& \mybullet \ar@{--}[ddd]^{>0} \ar@{--}[rrr]^{>0} 
\ar@{}[l]|-{\textstyle \tilde t\ } &&& 
\mybullet \ar@{--}[rrr]^{>0} &&&
\ar@{.}[rr] && \ar@{--}[rrr]^{>0} &&& 
\mybullet \ar@{--}[rrr]^{>0} &&& 
\mybullet \ar@{}[r]|-{\ \ \ \ \;\textstyle\tilde P_{n+1}} && \\
&&&&&&&&&&&&&&&&& \\
&&&&&&&&&&&&&&&&& \\
& \mybullet \ar@{}[r]|-{\ \ \ \textstyle\tilde Q_1} &&&&&&&&&&&&&&&& \\
&&&&&&&&&&&&&&&&& \\
}}
\end{equation}
Here $\tilde t$ is the unique vertex in $\tilde T$ where the shortest path from $\tilde Q_1$ to $\tilde P_2$ branches off towards $\tilde P_1$ or where it meets~$\tilde P_1$. The same properties characterize the vertex $t$ in the diagram (\ref{An1}); hence we must have $\tilde t=t$. By combining the diagrams (\ref{An1}) and (\ref{An2}) at this point it thus follows that $\tilde T$ contains a subtree of the form $A_{n+1}$, as desired.
\end{Proof}

%%%%%%%%%%%%%%%%%%%%%%%%%%%%%%%%%%%%%%%%%%%%%%

Now assume that $\tilde T$ contains a subtree of the form:
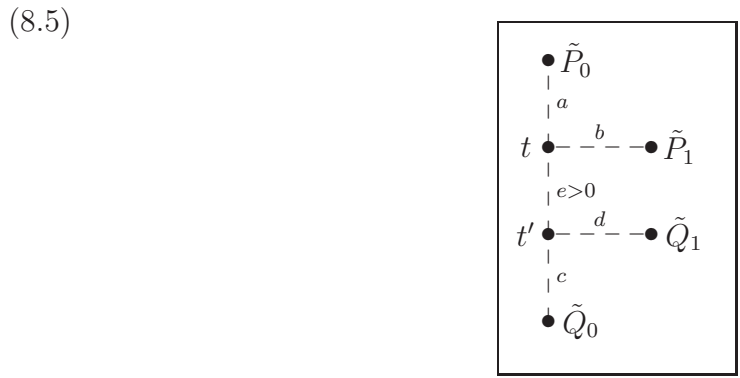
\UseTheoremCounterForNextEquation
\begin{equation}\label{bone1}
\fbox{\xymatrix@C=7pt@R=5pt{
&&&&&& \\
& \mybullet \ar@{--}[ddd]^{a} \ar@{}[r]|-{\ \ \textstyle\tilde P_0} &&&&& \\
&&&&&& \\
&&&&&& \\
& \mybullet \ar@{--}[ddd]^{e>0} \ar@{--}[rrr]^{b} 
\ar@{}[l]|-{\textstyle t\ } &&& 
\mybullet \ar@{}[r]|-{\ \ \textstyle\tilde P_1} && \\
&&&&&& \\
&&&&&& \\
& \mybullet \ar@{--}[ddd]^{c}  \ar@{--}[rrr]^{d} 
\ar@{}[l]|-{\textstyle t'\ } &&& 
\mybullet \ar@{}[r]|-{\ \ \ \textstyle\tilde Q_1} && \\
&&&&&& \\
&&&&&& \\
& \mybullet \ar@{}[r]|-{\ \ \ \textstyle\tilde Q_0} &&&&& \\
&&&&&& \\}}
\end{equation}
Our ultimate goal is to show that such a subtree is impossible, but that will require several steps. First observe that, as in the proof of Lemma \ref{An}, by Proposition \ref{SigmaTreeProp} (e) this subtree maps isomorphically to its image in $\tilde T/\gensigma$. In view of the right half of the diagram (\ref{SPlusSTildeT}) it follows that $T$ contains a subtree of the form:
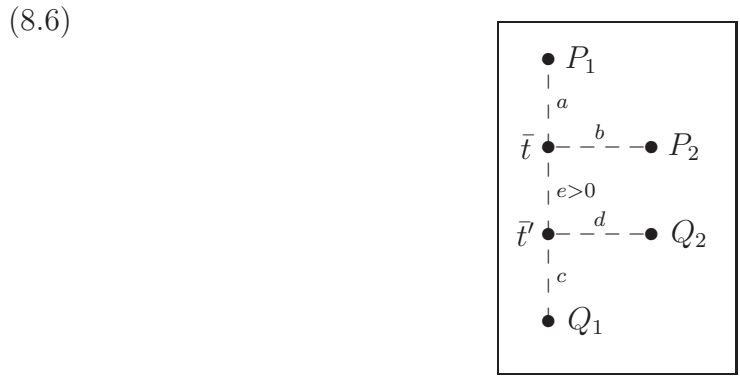
\UseTheoremCounterForNextEquation
\begin{equation}\label{bone2}
\fbox{\xymatrix@C=7pt@R=5pt{
&&&&&& \\
& \mybullet \ar@{--}[ddd]^{a} \ar@{}[r]|-{\ \ \ \textstyle P_1} &&&&& \\
&&&&&& \\
&&&&&& \\
& \mybullet \ar@{--}[ddd]^{e>0} \ar@{--}[rrr]^{b}
\ar@{}[l]|-{\textstyle \bar t\ } &&& 
\mybullet \ar@{}[r]|-{\ \ \ \textstyle P_2} && \\
&&&&&& \\
&&&&&& \\
& \mybullet \ar@{--}[ddd]^{c}  \ar@{--}[rrr]^{d}
\ar@{}[l]|-{\textstyle\bar t'\ } &&& 
\mybullet \ar@{}[r]|-{\ \ \ \ \textstyle Q_2} && \\
&&&&&& \\
&&&&&& \\
& \mybullet \ar@{}[r]|-{\ \ \ \ \textstyle Q_1} &&&&& \\
&&&&&& \\}}
\end{equation}
On the other hand $T$ can also be obtained from $\tilde T$ by stabilization with respect to the marking $\tilde s|\Gamma$. Here the vertices $t$ and $t'$ lie between $\tilde P_0$ and $\tilde Q_0$; hence we can apply Lemma \ref{SpineSurvives} to them. 

%%%%%%%%%%%%%%%%%%%%%%%%%%%%%%%%%%%%%%%%%%%%%%

\medskip
{\bf Case 1: The vertex $t$ survives in~$T$:} In this case, constructing $T$ from $\tilde T$ by stabilization, we deduce from (\ref{bone1}) that $T$ contains a subtree of the form
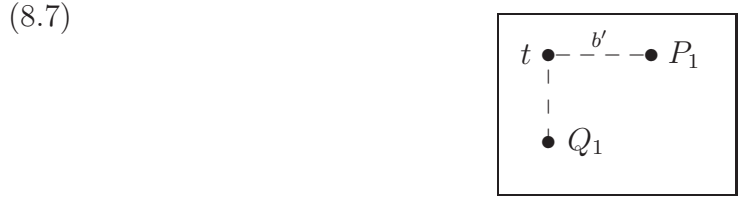
\UseTheoremCounterForNextEquation
\begin{equation}\label{bone3}
\fbox{\xymatrix@C=7pt@R=5pt{
&&&&&& \\
& \mybullet \ar@{--}[ddd] \ar@{--}[rrr]^{b'} 
\ar@{}[l]|-{\textstyle t\ } &&& 
\mybullet \ar@{}[r]|-{\ \ \ \textstyle P_1} && \\
&&&&&& \\
&&&&&& \\
& \mybullet \ar@{}[r]|-{\ \ \ \ \textstyle Q_1} &&&&& \\
&&&&&& \\}}
\end{equation}
with $0\le b'\le b$.

\begin{Lem}\label{B'GreaterA}
The case that $t$ survives in $T$ and that $b'>a$ is impossible.
\end{Lem}

\begin{Proof}
Suppose that $b'>a$. Then $b\ge b'>a\ge0$ implies that $b>0$. Moving along the shortest path in $T$ from $P_1$ towards~$Q_1$ in the diagrams (\ref{bone2}) and (\ref{bone3}), the inequality $b'>a$ shows that we first reach the branch off point towards $P_2$ before reaching~$t$. This implies that $T$ contains a subtree of the form:
$$\fbox
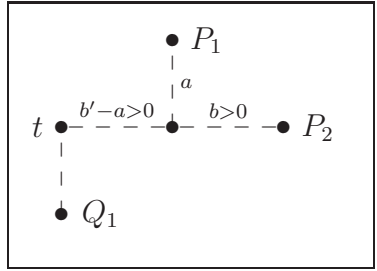
{\xymatrix@C=8pt@R=5pt{
&&&&&&&& \\
&&&& \mybullet \ar@{--}[ddd]^{a} \ar@{}[r]|-{\ \ \ \textstyle P_1} &&&& \\
&&&&&&&& \\
&&&&&&&& \\
& \mybullet \ar@{--}[ddd] \ar@{--}[rrr]^{b'-a>0}
\ar@{}[l]|-{\textstyle t\ } &&& 
\mybullet \ar@{--}[rrr]^{b>0} &&&
\mybullet \ar@{}[r]|-{\ \ \ \textstyle P_2} & \\
&&&&&&&& \\
&&&&&&&& \\
& \mybullet \ar@{}[r]|-{\ \ \ \ \textstyle Q_1} &&&&&&& \\
&&&&&&&&& \\}}$$
Since $\tilde T$ can be obtained from $T$ by inserting edges, it follows that $\tilde T$ contains a subtree of the form:
$$\fbox
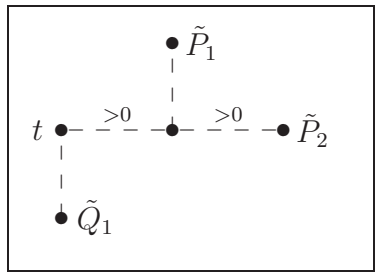
{\xymatrix@C=8pt@R=5pt{
&&&&&&&& \\
&&&& \mybullet \ar@{--}[ddd] \ar@{}[r]|-{\ \ \textstyle\tilde P_1} &&&& \\
&&&&&&&& \\
&&&&&&&& \\
& \mybullet \ar@{--}[ddd] \ar@{--}[rrr]^{>0}
\ar@{}[l]|-{\textstyle t\ } &&& 
\mybullet \ar@{--}[rrr]^{>0} &&&
\mybullet \ar@{}[r]|-{\ \ \textstyle\tilde P_2} & \\
&&&&&&&& \\
&&&&&&&& \\
& \mybullet \ar@{}[r]|-{\ \ \ \textstyle\tilde Q_1} &&&&&&& \\ 
&&&&&&&&& \\}}$$
Combining this with the subtree (\ref{bone1}) it follows that $\tilde T$ contains a subtree of the form forbidden by Lemma \ref{An}. Thus the case is impossible.
\end{Proof}

\begin{Lem}\label{B'EqualA}
The case that both $t$ and $t'$ survive in $T$ and that $b'=a$ is impossible.
\end{Lem}

\begin{Proof}
Suppose that $b'=a$. Measuring the distances along the shortest path in $T$ from $P_1$ towards~$Q_1$ in the diagrams (\ref{bone2}) and (\ref{bone3}), we then find that $\bar t=t$. Let $t''$ be the vertex in $\tilde T$ adjacent to $t$ in the direction towards~$t'$. Then either $t''=t'$ and it survives in $T$ by assumption, or it lies strictly between $t$ and~$t'$, and hence also strictly between $\tilde P_0$ and~$\tilde Q_0$, in which case it survives in $T$ by Lemma \ref{SpineSurvives} (a). By comparing the diagrams  (\ref{bone1}) and (\ref{bone2}) we then find that $t''$ is mapped to itself under the projection $\tilde T\onto \tilde T/\gensigma\cong T$. It follows that the edge $(t,t'')$ of $\tilde T$ survives in $T$ under stabilization and is mapped to itself under the projection $\tilde T\onto \tilde T/\gensigma\cong T$.

By the definition of the dual tree, the edge $(t,t'')$ represents a singular point $\tilde x_0$ of the closed fiber of~$\tilde X$, respectively a singular point $x_0$ of the closed fiber of~$X$. The fact that the edge survives under stabilization means that the stabilization morphism $\tilde X\to X$ maps $\tilde x_0$ to $x_0$ and is a local isomorphism there. The fact that the edge is mapped to itself under the projection $\tilde T\onto \tilde T/\gensigma\cong T$ means that the morphism $\tilde X\onto\tilde X/\gensigma\cong X$ also maps $\tilde x_0$ to~$x_0$.
By Proposition \ref{QuotientProp} (b) the latter implies that $\tilde X$ is \'etale locally isomorphic to $\Spec R[y,z]/(yz-\pi^n)$ at~$\tilde x_0$, and that $X$ is \'etale locally isomorphic to $\Spec R[u,v]/(uv-\pi^{2n})$ at~$x_0$, for some integer $n\ge1$. As the other morphism $\tilde X\to X$ is a local isomorphism at~$\tilde x_0$, it follows that the two charts are \'etale locally isomorphic at the singular point of the fiber. But the local equations $yz=\pi^n$ and $uv=\pi^{2n}$ are not equivalent for $n\ge1$; hence the case is not possible.
\end{Proof}

\begin{Lem}\label{BothSurvive}
The case that both $t$ and $t'$ survive in $T$ is impossible.
\end{Lem}

\begin{Proof}
Suppose that $t$ and $t'$ both survive in~$T$. Then any vertex strictly between $t$ and $t'$ lies strictly between $\tilde P_0$ and $\tilde Q_0$ and hence also survives in~$T$ by Lemma \ref{SpineSurvives} (a). Thus from (\ref{bone1}) we deduce that $T$ contains a subtree of the form
\UseTheoremCounterForNextEquation
\begin{equation}\label{bone3x}
\fbox{\xymatrix@C=7pt@R=5pt{
&&&&&& \\
& \mybullet \ar@{--}[ddd]^{e>0} \ar@{--}[rrr]^{b'} 
\ar@{}[l]|-{\textstyle t\ } &&& 
\mybullet \ar@{}[r]|-{\ \ \ \textstyle P_1} && \\
&&&&&& \\
&&&&&& \\
& \mybullet \ar@{--}[rrr]^{d'} 
\ar@{}[l]|-{\textstyle t'\ } &&& 
\mybullet \ar@{}[r]|-{\ \ \ \ \textstyle Q_1} && \\
&&&&&& \\}}
\end{equation}
with $0\le b'\le b$ and $0\le d'\le d$, but where where the distance $e$ between $t$ and $t'$ has not changed. Comparing the distances between $P_1$ and $Q_1$ in both (\ref{bone2}) and (\ref{bone3x}) we find that $a+e+c = b'+e+d'$ and hence  
\UseTheoremCounterForNextEquation
\begin{equation}\label{myeq1}
a+c = b'+d'.
\end{equation}
By Lemmas \ref{B'GreaterA} and \ref{B'EqualA} we already know that $b'<a$. On the other hand, the present situation is invariant under interchanging the two critical points $\tilde s(i_0)$ and $\tilde s(j_0)$. This interchanges $\tilde P_i$ with~$\tilde Q_i$ and flips the diagrams (\ref{bone1})  and (\ref{bone3x}) upside down. In particular, it interchanges the vertices $t$ and~$t'$ and the distances $(a,b')$ and $(c,d')$. Thus by symmetry, the proof of $b'<a$ also shows that $d'<c$. But together these inequalities imply that $b'+d' < a+c$, which contradicts the equality (\ref{myeq1}). Thus the case that both $t$ and $t'$ survive in $T$ is impossible.
\end{Proof}

%%%%%%%%%%%%%%%%%%%%%%%%%%%%%%%%%%%%%%%%%%%%%%

\medskip
{\bf Case 2: The vertex $t$ does not survive in~$T$:} 
By Lemma \ref{SpineSurvives} this requires that $t=\tilde P_0$ and $a=0$, so that the subtree of $\tilde T$ in (\ref{bone1}) really looks like this:
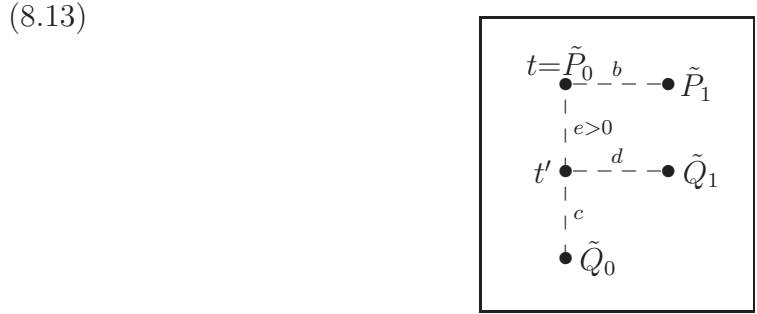
\UseTheoremCounterForNextEquation
\begin{equation}\label{bone1a}
\fbox{\xymatrix@C=7pt@R=5pt{
&&&&&&& \\
&&&&&&& \\
&& \mybullet \ar@{--}[ddd]^{e>0} \ar@{--}[rrr]^{b} 
\ar@{}[uu]|-{\textstyle t{=}\tilde P_0\ } &&& 
\mybullet \ar@{}[r]|-{\ \ \textstyle\tilde P_1} && \\
&&&&&&& \\
&&&&&&& \\
&& \mybullet \ar@{--}[ddd]^{c}  \ar@{--}[rrr]^{d} 
\ar@{}[l]|-{\textstyle t'\ } &&& 
\mybullet \ar@{}[r]|-{\ \ \ \textstyle\tilde Q_1} && \\
&&&&&&& \\
&&&&&&& \\
&& \mybullet \ar@{}[r]|-{\ \ \ \textstyle\tilde Q_0} &&&&& \\
&&&&&&& \\}}
\end{equation}
Then the subtree of $T$ in (\ref{bone2}) looks like this:
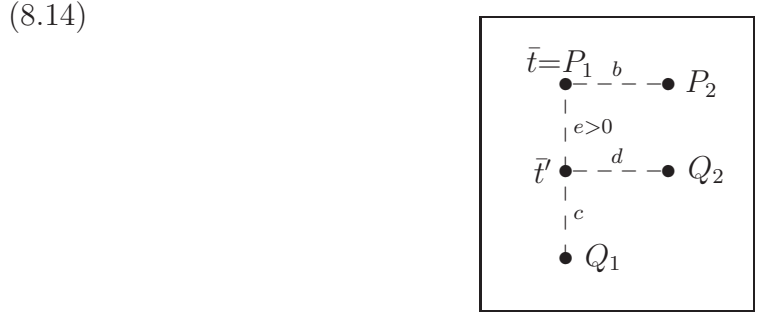
\UseTheoremCounterForNextEquation
\begin{equation}\label{bone2a}
\fbox{\xymatrix@C=7pt@R=5pt{
&&&&&&& \\
&&&&&&& \\
&& \mybullet \ar@{--}[ddd]^{e>0} \ar@{--}[rrr]^{b} 
\ar@{}[uu]|-{\textstyle \bar t{=}P_1\ } &&& 
\mybullet \ar@{}[r]|-{\ \ \ \textstyle P_2} && \\
&&&&&&& \\
&&&&&&& \\
&& \mybullet \ar@{--}[ddd]^{c}  \ar@{--}[rrr]^{d}
\ar@{}[l]|-{\textstyle\bar t'\ } &&& 
\mybullet \ar@{}[r]|-{\ \ \ \ \textstyle Q_2} && \\
&&&&&&& \\
&&&&&&& \\
&& \mybullet \ar@{}[r]|-{\ \ \ \ \textstyle Q_1} &&&&& \\
&&&&&& \\}}
\end{equation}

\begin{Lem}\label{BGreater0}
The case that $t$ does not survive in $T$ and that $b>0$ is impossible.
\end{Lem}

\begin{Proof}
Suppose that $b>0$. Since $t=\tilde P_0$ does not survive, we must then be in the case (ii) of Lemma \ref{SpineSurvives} (b), with $\tilde P_1$ in the connected component~$S$ of $\tilde T\setminus\{\tilde P_0\}$. The condition in \ref{SpineSurvives} (b) (ii) implies that all markings of vertices in $S$ are of the form $\tilde s(\gamma)$ for $\gamma\in\Gamma$. This implies that all vertices in $S$ satisfy the condition in Proposition \ref{StabilizationProp} (a) and therefore survive in~$T$. Thus on passing from $\tilde T$ to~$T$, the whole subtree $S$ is preserved. In particular $\tilde P_1=P_1$ survives in~$T$, and all distances within $S$ are preserved. 

Let $t'''$ be the vertex in $S$ connected by an edge to $t=\tilde P_0$. Since $t'''$ survives but $t$ does not, the edge $(t''',t)$ is either moved from $t$ to another vertex $\tilde t$ in $\tilde T\setminus S$ or removed entirely. 
Suppose that during the stabilization process, some marking in $\tilde T\setminus S$ is moved into~$S$. Since marked points are moved only when end vertices, i.e., leaves of the tree are removed, this requires that the whole complement $\tilde T\setminus S$ is pushed onto~$t'''$. Thus we have two possibilities: 
Either the edge $(t''',t)$ is replaced by an edge $(t''',\tilde t\,)$ for some $\tilde t$ outside~$S$ and all markings of $\tilde T$ outside $S$ remain outside~$S$.
Or the edge $(t''',t)$ is removed and all markings of $\tilde T$ outside $S$ are moved onto~$t'''$.

In either case we observe that $t'''$ lies between $Q_1$ and~$P_1$. Note also that the diagram (\ref{bone2a}) and the assumption $b>0$ show that $P_1$ lies strictly between $Q_1$ and $P_2$ in~$T$. Thus we have the following subtree of $T$:
\UseTheoremCounterForNextEquation
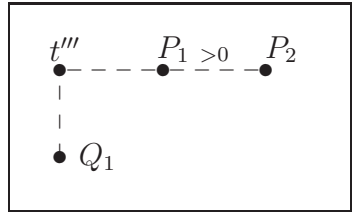
\begin{equation}\label{An2tutu}
\fbox{\xymatrix@C=7pt@R=5pt{
&&&&&&&&& \\
&&&&&&&&& \\
& \mybullet \ar@{--}[ddd] \ar@{--}[rrr] 
\ar@{}[uu]|-{\ \textstyle t'''}
&&& \mybullet \ar@{--}[rrr]^{>0} 
\ar@{}[uu]|-{\ \ \textstyle P_1} &&& 
\mybullet 
% \ar@{}[r]|-{\ \ \ \textstyle P_2} 
\ar@{}[uu]|-{\ \ \ \textstyle P_2}
&& \\
&&&&&&&&& \\
&&&&&&&&& \\
& \mybullet \ar@{}[r]|-{\ \ \ \ \textstyle Q_1} &&&&&&&& \\
&&&&&&&&& \\}}
\end{equation}
This shows that $P_2$ lies in~$S$ and that $Q_1$, $P_1$, $P_2$ are all distinct.

If $\tilde P_2$ does not lie in~$S$, we must have the second case of the above remarks, and so during stabilization, the points $\tilde P_2$ and $\tilde Q_1$ are both moved to~$t'''$. This contradicts the fact that $P_2\not=Q_1$. Thus $\tilde P_2$ lies in~$S$ and is therefore equal to~$P_2$.

Moreover, the subtree (\ref{An2tutu}) shows that $P_1$ lies between $t'''$ and $P_2$ within the image of~$S$. As $S$ is mapped isomorphically to its image, it follows that $\tilde P_1$ lies between $t'''$ and $\tilde P_2$ within~$S$. Thus $\tilde T$ contains a subtree of the form:
$$\fbox
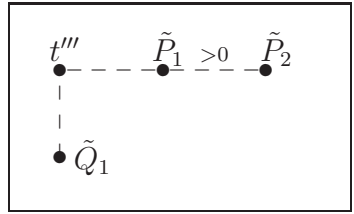
{\xymatrix@C=7pt@R=5pt{
&&&&&&&&& \\
&&&&&&&&& \\
& \mybullet \ar@{--}[ddd] \ar@{--}[rrr] 
\ar@{}[uu]|-{\ \textstyle t'''}
&&& \mybullet \ar@{--}[rrr]^{>0} 
\ar@{}[uu]|-{\ \textstyle\tilde P_1} &&& 
\mybullet 
% \ar@{}[r]|-{\ \ \ \textstyle P_2} 
\ar@{}[uu]|-{\ \ \textstyle\tilde P_2}
&& \\
&&&&&&&&& \\
&&&&&&&&& \\
& \mybullet \ar@{}[r]|-{\ \ \ \textstyle\tilde Q_1} &&&&&&&& \\
&&&&&&&&& \\}}$$
Combining this with the diagram (\ref{bone1a}) we deduce that $\tilde T$ contains a subtree of the form:
$$\fbox
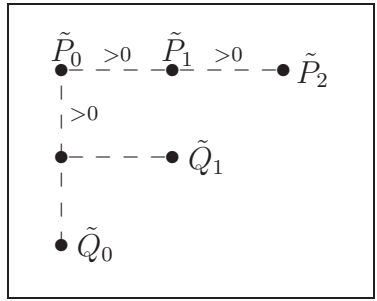
{\xymatrix@C=8pt@R=5pt{
&&&&&&&& \\
&&&&&&&& \\
& \mybullet \ar@{--}[ddd]^{>0} \ar@{--}[rrr]^{>0} 
\ar@{}[uu]|-{\ \textstyle\tilde P_0} &&& 
\mybullet \ar@{--}[rrr]^{>0} \ar@{}[uu]|-{\ \textstyle\tilde P_1} &&&
\mybullet \ar@{}[r]|-{\ \ \textstyle\tilde P_2} & \\
&&&&&&&& \\
&&&&&&&& \\
& \mybullet \ar@{--}[ddd]  \ar@{--}[rrr] &&& 
\mybullet \ar@{}[r]|-{\ \ \ \textstyle\tilde Q_1} &&&& \\
&&&&&&&& \\
&&&&&&&& \\
& \mybullet \ar@{}[r]|-{\ \ \ \textstyle\tilde Q_0} &&&&&&& \\
&&&&&&&&& \\}}$$
But this is a special case of the subtree excluded by Lemma \ref{An}. Thus the case under consideration is impossible.
\end{Proof}

%%%%%%%%%%%%%%%%%%%%%%%%%%%%%%%%%%%%%%%%%%%%%%

\begin{Lem}\label{OneSurvives}
The case that precisely one of $t$, $t'$ survives in $T$ is impossible.
\end{Lem}

\begin{Proof}
By the same symmetry as in the proof of Lemma \ref{BothSurvive}, we may without loss of generality assume that $t$ survives in~$T$ while $t'$ does not. 
Then by Lemma \ref{B'GreaterA} we must have $b'\le a$. On the other hand, by Lemma \ref{BGreater0} and the remarks preceding it, and by symmetry, we must also have $c=d=0$. Thus the subtree of $\tilde T$ in the diagram (\ref{bone1}) looks like this:
$$\fbox
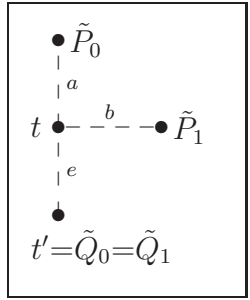
{\xymatrix@C=7pt@R=5pt{
&&&&&& \\
& \mybullet \ar@{--}[ddd]^{a} \ar@{}[r]|-{\ \ \textstyle\tilde P_0} &&&&& \\
&&&&&& \\
&&&&&& \\
& \mybullet \ar@{--}[ddd]^{e} \ar@{--}[rrr]^{b} 
\ar@{}[l]|-{\textstyle t\ } &&& 
\mybullet \ar@{}[r]|-{\ \ \textstyle\tilde P_1} && \\
&&&&&& \\
&&&&&& \\
& \mybullet
\ar@{}[ddrr]|-{\ \textstyle t'{=}\tilde Q_0{=}\tilde Q_1} &&&&& \\
&&&&&& \\
&&&&&& \\}}$$
Consequently, the subtree of $T$ obtained by taking quotients in (\ref{bone2}) looks like this:
$$\fbox
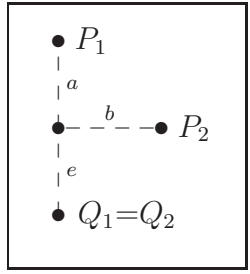
{\xymatrix@C=7pt@R=5pt{
&&&&&& \\
& \mybullet \ar@{--}[ddd]^{a} \ar@{}[r]|-{\ \ \ \textstyle P_1} &&&&& \\
&&&&&& \\
&&&&&& \\
& \mybullet \ar@{--}[ddd]^{e} \ar@{--}[rrr]^{b} &&& 
\mybullet \ar@{}[r]|-{\ \ \ \textstyle P_2} && \\
&&&&&& \\
&&&&&& \\
& \mybullet 
\ar@{}[rrrr]|-{\textstyle Q_1{=}Q_2} &&&&& \\
&&&&&& \\}}$$
Moreover, the subtree of $T$ obtained by stabilization in (\ref{bone3}) looks like this:
$$\fbox
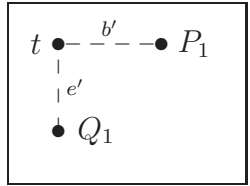
{\xymatrix@C=7pt@R=5pt{
&&&&&& \\
& \mybullet \ar@{--}[ddd]^{e'} \ar@{--}[rrr]^{b'} 
\ar@{}[l]|-{\textstyle t\ } &&& 
\mybullet \ar@{}[r]|-{\ \ \ \textstyle P_1} && \\
&&&&&& \\
&&&&&& \\
& \mybullet \ar@{}[r]|-{\ \ \ \ \textstyle Q_1} &&&&& \\
&&&&&& \\}}$$
Here $0\le e'<e$, because the vertex $t'=\tilde Q_1$ does not survive in~$T$.
Comparing the distances between $P_1$ and $Q_1$ in the last two diagrams we find that $a+e = b'+e'$. But since $b'\le a$ and $e'< e$, we also have $b'+e'<a+e$. Together this gives a contradiction; hence the case under consideration is not possible.
\end{Proof}

%%%%%%%%%%%%%%%%%%%%%%%%%%%%%%%%%%%%%%%%%%%%%%

\begin{Lem}\label{NoneSurvives}
The case that none of $t$, $t'$ survives in $T$ is impossible.
\end{Lem}

\begin{Proof}
If none of $t$, $t'$ survives, by Lemma \ref{BGreater0} and the remarks preceding it, and by symmetry, we must have $a=b=c=d=0$. Thus the subtree of $\tilde T$ in the diagram (\ref{bone1}) looks like this:
$$\fbox
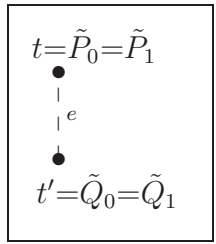
{\xymatrix@C=7pt@R=5pt{
&&&&& \\
&&&&& \\
& \mybullet \ar@{--}[ddd]^{e}
\ar@{}[uurr]|-{\ \textstyle t{=}\tilde P_0{=}\tilde P_1} &&&& \\
&&&&& \\
&&&&& \\
& \mybullet
\ar@{}[ddrr]|-{\ \ \textstyle t'{=}\tilde Q_0{=}\tilde Q_1} &&&& \\
&&&&& \\
&&&&& \\}}$$
Consequently, the subtree of $T$ obtained by taking quotients in (\ref{bone2}) looks like this:
$$\fbox
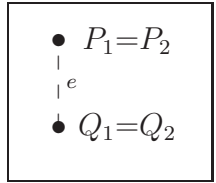
{\xymatrix@C=7pt@R=5pt{
&&&&& \\
& \mybullet \ar@{--}[ddd]^{e} 
\ar@{}[rrrr]|-{\textstyle P_1{=}P_2} &&&& \\
&&&&& \\
&&&&& \\
& \mybullet 
\ar@{}[rrrr]|-{\textstyle Q_1{=}Q_2} &&&& \\
&&&&& \\}}$$
Moreover, the subtree of $T$ obtained by stabilization has the form
$$\fbox
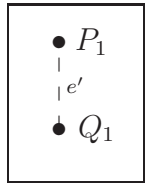
{\xymatrix@C=7pt@R=5pt{
&&& \\
& \mybullet \ar@{--}[ddd]^{e'} \ar@{}[r]|-{\ \ \ \textstyle P_1} && \\
&&& \\
&&& \\
& \mybullet \ar@{}[r]|-{\ \ \ \ \textstyle Q_1} && \\
&&& \\}}$$
Here $0\le e'\le e-2$, because both $t$ and $t'$ do not survive in~$T$. But comparing the distances between $P_1$ and $Q_1$ in the last two diagrams we find that $e=e'$, which yields a contradiction. Thus the case under consideration is not possible.
\end{Proof}

By combining Lemmas \ref{BothSurvive}, \ref{OneSurvives}, and \ref{NoneSurvives} we thus arrive at the following conclusion:

\begin{Thm}\label{NoSeparation}
The tree $\tilde T$ does not contain a subtree of the form (\ref{bone1}):
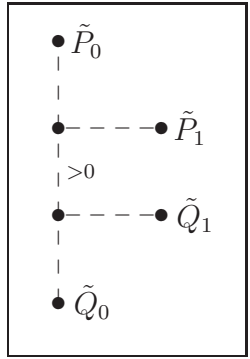
$$\fbox{\xymatrix@C=7pt@R=5pt{
&&&&&& \\
& \mybullet \ar@{--}[ddd] \ar@{}[r]|-{\ \ \textstyle\tilde P_0} &&&&& \\
&&&&&& \\
&&&&&& \\
& \mybullet \ar@{--}[ddd]^{>0} \ar@{--}[rrr] &&& 
\mybullet \ar@{}[r]|-{\ \ \textstyle\tilde P_1} && \\
&&&&&& \\
&&&&&& \\
& \mybullet \ar@{--}[ddd]  \ar@{--}[rrr] &&& 
\mybullet \ar@{}[r]|-{\ \ \ \textstyle\tilde Q_1} && \\
&&&&&& \\
&&&&&& \\
& \mybullet \ar@{}[r]|-{\ \ \ \textstyle\tilde Q_0} &&&&& \\
&&&&&& \\}}$$
\end{Thm}

%%%%%%%%%%%%%%%%%%%%%%%%%%%%%%%%%%%%%%%%%%%%%%%%%%%%%%%%%%%%%%%%%%%%%%%%%%%%%%%%%%%%%%%%%%%%%

\section{Proof of Theorem \ref{GammaModuliQuasiFinite}}
\label{ProofOfTheorem}

We must show that the fiber of the moduli space $M_\Gamma$ over any field $k$ of characteristic $\not=2$ is finite. Suppose that this is not the case. Then the fiber contains an irreducible affine curve~$D$. After extending $k$ and shrinking~$D$, if necessary, we may assume that $D$ is smooth over~$k$. Let $(C,f,P,Q,s)$ be the pullback to $D$ of the universal $\Gamma$-marked quadratic morphism over~$M_\Gamma$. Let $\bar D$ be a smooth compactification of $D$ over~$k$. Then over the local ring at each point of $\bar D\setminus D$ we can study the degeneration as in Sections \ref{StableCurvesQuad} and \ref{Exclude}. The key consequence is this:
 
\begin{Lem}\label{CrossRatio}
The cross ratio of $P$, $f(P)$, $Q$, $f(Q)$ is constant over~$D$.
\end{Lem}

\begin{Proof}
The cross ratio of four points in $\BP^1$ induces an isomorphism from the fine moduli space of stable $4$-pointed curves of genus zero to~$\BP^1$. Here the smooth curves correspond to cross ratios $\not=0$, $1$, $\infty$, while the three exceptional values correspond to degenerate curves consisting of two rational curves connected by an ordinary double point, each of which carries two marked points. The three degenerate curves correspond to the three ways of separating the four points into two pairs.

In our case, if the elements $i_0$, $i_1$, $j_0$, $j_1$ of $\tilde\Gamma$ are not all distinct, some of the corresponding marked points are equal, and so their cross ratio is $0$, $1$, or~$\infty$, in which case we are done. So assume that $i_0$, $i_1$, $j_0$, $j_1$ are distinct. 

Then the corresponding sections are disjoint, and so their cross ratio defines a morphism $D\to\BP^1_k\setminus\{0,1,\infty\}$, which classifies the associated smooth stable curve with the four distinct marked points $P$, $f(P)$, $Q$, $f(Q)$. Assume that this morphism is not constant. Then the unique extension $\bar D\to\BP^1_k$ is surjective. Pulling back the universal family from the fine moduli space $\BP^1_k$ we obtain a stable $4$-pointed curve of genus zero over~$\bar D$. By surjectivity, all possible degenerations must occur. In particular, there is a point in $\bar D\setminus D$ where the points $P$ and $f(P)$ lie in one irreducible component of the fiber, while $Q$ and $f(Q)$ lie in the other. Let $R$ denote the local ring of $\bar D$ at this point and $X'$ the pullback of the universal family to~$\Spec R$ and $s'$ the associated marking indexed by $\{i_0,i_1,j_0,j_1\} \subset \tilde\Gamma$. Then the dual tree $T'$ of the closed fiber of~$X'$ looks like this:
$$\fbox
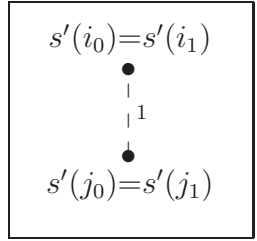
{\xymatrix@C=7pt@R=5pt{
&&&& \\
&&&& \\
&&& \mybullet \ar@{--}[ddd]^{1}
\ar@{}[uur]|{\textstyle s'(i_0){=}s'(i_1)\quad } & \\
&&&& \\
&&&& \\
&&& \mybullet
\ar@{}[ddr]|{\textstyle s'(j_0){=}s'(j_1)\quad } & \\
&&&& \\
&&&& \\}}$$

Now let $(\tilde X,\tilde s)$ be the stable marked curve over $\Spec R$ described in Section \ref{StableCurvesQuad}. By construction, the generic fiber of $(X',s')$ is obtained from the generic fiber of $(\tilde X,\tilde s)$ by forgetting some marked points. Thus the whole family $(X',s')$ is obtained from $(\tilde X,\tilde s)$ by the stabilization process described in Section \ref{StableCurves}. In the same way the dual tree $T'$ is obtained from the dual tree~$\tilde T$.

The above form of $T'$ shows that the markings $\tilde P_0 = \tilde s(i_0)$ and $\tilde P_1 = \tilde s(i_1)$ of $\tilde T$ are moved to one vertex of~$T'$, and the markings $\tilde Q_0 = \tilde s(j_0)$ and $\tilde Q_1 = \tilde s(j_1))$ are moved to a different vertex of~$T'$. Since we can, conversely, obtain $\tilde T$ from $T'$ by inserting edges and thereby moving certain marked points apart or further apart, but never the opposite, it follows that $\tilde T$ contains a subtree of the form:
$$\fbox
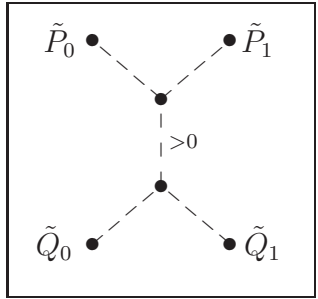
{\xymatrix@C=7pt@R=5pt{
&&&&&&&& \\
&& \mybullet \ar@{--}[ddrr] \ar@{}[l]|-{\textstyle\tilde P_0\ \ \ } 
&&&& \mybullet \ar@{}[r]|-{\ \ \textstyle\tilde P_1} && \\
&&&&&&&& \\
&&&& \mybullet \ar@{--}[ddd]^{>0} \ar@{--}[rruu] &&&& \\
&&&&&&&& \\
&&&&&&&& \\
&&&& \mybullet \ar@{--}[rrdd] &&&& \\
&&&&&&&& \\
&& \mybullet \ar@{--}[uurr] \ar@{}[l]|-{\textstyle\tilde Q_0\ \ \ \ } 
&&&& \mybullet \ar@{}[r]|-{\ \ \ \textstyle\tilde Q_1} && \\
&&&&&&&& \\}}$$
But this contradicts Theorem \ref{NoSeparation}.
\end{Proof}

%%%%%%%%%%%%%%%%%%%%%%%%%%%%%%%%%%%%%%%%%%%%%%

\medskip
The remaining arguments are relatively pedestrian. (Although I expect that with the methods of the preceding section one can also exclude certain subtrees of $\tilde T$ containing the marked points $\tilde P_0$, $\tilde P_1$, $\tilde P_2$, $\tilde Q_0$, $\tilde Q_1$, and from this deduce that the parameter $t$ discussed below is constant.)

\medskip
Since $|\Gamma|\ge3$ by assumption, using Proposition \ref{Unstable} we may without loss of generality assume that $i_1\not=i_0,j_0$. Then as in the proof of Proposition \ref{StableModuli} we may assume that $C=\BP^1\times D$ with $\tilde s(i_0)=0$ and $\tilde s(j_0)=\infty$ and $\tilde s(i_1)=1$, and hence $f(x) = \frac{cx^2+1}{dx^2+1}$ for $c$, $d\in\Gamma(D,\CO_D)$. Then $\tilde s(j_1)=f(\infty)= \frac{c}{d}$, and by Lemma \ref{CrossRatio} this is a constant section $C\to\BP^1$. Since $f$ is a well-defined quadratic morphism in every fiber, we also have $c\not=d$ everywhere. Thus in projective coordinates we can write $(c:d) = (\alpha:\beta)$ for constants $\alpha\not=\beta$. This means that we can substitute $c=\alpha t$ and $d=\beta t$ for some function $t$ on~$C$. Thereafter we have
\UseTheoremCounterForNextEquation
\begin{equation}\label{FOneParameter}
f(x) = \frac{\alpha tx^2+1}{\beta tx^2+1}.
\end{equation}

\begin{Lem}\label{FinEq}
There exist $k>\ell>0$ with $f^k(0)=-f^\ell(0)$.
\end{Lem}

\begin{Proof}
As the postcritical orbit of $f$ is finite, there exist $n>m\ge0$ with $f^n(0)=f^m(0)$. Among these pairs $(n,m)$ select one where $m$ is minimal. 
If $m=0$, the equation simply reads $f^n(0)=0$ and implies that $f^{2n}(0) = 0 = -f^n(0)$, and so the desired assertion holds with $(k,\ell)=(2n,n)$. 
If $m=1$, the equation reads $f^n(0)=f(0)$, which by (\ref{RamLift}) implies that $f^{n-1}(0)=0$, contradicting the minimality of~$m$.
So assume that $m\ge2$. Then by minimality we have $f^n(0)=f^m(0)$ and $f^{n-1}(0)\not=f^{m-1}(0)$. Since $f$ is a Galois covering of degree $2$ with the non-trivial Galois automorphism $x\mapsto-x$, this implies that $f^{n-1}(0)=-f^{m-1}(0)$. Thus the desired assertion holds with $(k,\ell)=(n-1,m-1)$.
\end{Proof}

\medskip
To finish the proof of Theorem \ref{GammaModuliQuasiFinite}, we must show that for any $\Gamma$-marked quadratic morphism of the form (\ref{FOneParameter}) over a curve over a field, where $\alpha\not=\beta$ are fixed constants, the parameter $t$ is constant as well. We will achieve this by writing out the equation $f^k(0)=-f^\ell(0)$ in terms of $t$ and showing that it has only finitely many solutions over any field of characteristic $\not=2$.

\medskip
{\bf Case 1:} $\beta=0$.\ \ 
This is the case where $f(\infty)=\infty$, that is, where $j_0=j_1\in\Gamma$. Here we may without loss of generality assume that $\alpha=1$, so that $f(x) = tx^2+1$. Then $g_n := f^n(0) \in \BZ[t]$ is characterized recursively by
$$g_0 :=0 \qquad\hbox{and}\qquad
  g_{n+1} := tg_n^2+1 \quad \hbox{for all $n\ge0$.}$$
By an easy induction this implies that for every $n\ge1$ the polynomial $g_n$ is monic of degree $2^{n-1}-1$ in~$t$. The equation $f^k(0)=-f^\ell(0)$ from Lemma \ref{FinEq} now has the form $g_k+g_\ell=0$. Since $k>\ell>0$, this equation is monic of degree $2^{k-1}-1$ in~$t$. It therefore has only finitely many solutions over any field, as desired.

\medskip
{\bf Case 2:} $\alpha=0$.\ \ 
This is the case where $f(\infty)=0$, that is, where $i_0=j_1\in\Gamma$. Here we may without loss of generality assume that $\beta=1$, so that $f(x) = \frac{1}{tx^2+1}$. In projective coordinates we can then write $f^n(0)=(g_n:h_n)$ where $g_n$, $h_n\in\BZ[t]$ are characterized recursively by
$$\left\{\begin{array}{l}
g_0:=0 \\[6pt]
h_0:=1
\end{array}\right\}
\quad\hbox{and}\quad
\left\{\begin{array}{l}
g_{n+1} := h_n^2 \\[6pt]
h_{n+1} := tg_n^2+h_n^2
\end{array}\right\}
\quad \hbox{for all $n\ge0$.}$$

\begin{Lem}\label{MonicCase2}
Here for every $n\ge0$ the polynomial $h_n$ is monic in~$t$.
\end{Lem}

% Actually the degree of $h_n$ is $\lfloor\frac{2^n}{3}\rfloor$.

\begin{Proof}
The recursion formula implies that $h_0=h_1=1$ and $h_{n+2} = th_n^4+h_{n+1}^2$ for all $n\ge0$. If $h_n$ and $h_{n+1}$ are monic, then $th_n^4$ and $h_{n+1}^2$ are monic of odd, respectively, even degree. Then $h_{n+2}$ is monic of the greater of these two degrees, and the lemma follows by induction.
\end{Proof}

In projective coordinates the equation $f^k(0)=-f^\ell(0)$ now reads ${(g_k:h_k)} = {(-g_\ell:h_\ell)}$, which means that $g_k h_\ell + g_\ell h_k=0$. By the recursion formula this is equivalent to 
$$h_{k-1}^2 h_\ell + h_{\ell-1}^2 h_k=0.$$
Lemma \ref{MonicCase2} implies that both summands are monic in $t$ of certain degrees. If these degrees are equal, the highest term in $t$ of the equation has the form~$2t^N$, while otherwise it has the form~$t^N$, for some $N\ge0$. 
% And both of these cases occur, depending on the parities of $k$ and $\ell$.
In either case the equation has only finitely many solutions over any field of characteristic $\not=2$, as desired.

\medskip
{\bf Case 3:} $\alpha$, $\beta\not=0$.\ \ 
This is the case where $f(\infty)=\frac{\alpha}{\beta} \not=0$, $\infty$. Since $\alpha\not=\beta$ we also have $f(\infty)\not=1$, so that $\tilde s(i_0)=0$ and $\tilde s(i_1)=1$ and $\tilde s(j_0)=\infty$ and $\tilde s(j_1)=f(\infty)=\frac{\alpha}{\beta}$ are all distinct.
% (Thus $i_0$, $i_1$, $j_0$, $j_1$ of $\tilde\Gamma$ are all distinct. In any case the equality $i_1=j_1$ is already forbidden by the condition \ref{AbsPostCritDef} (c).)
In projective coordinates we can write $f^n(0)=(g_n:h_n)$ where $g_n$, $h_n\in\BZ[\alpha,\beta,t]$ are characterized recursively by
$$\left\{\begin{array}{l}
g_0:=0 \\[6pt]
h_0:=1
\end{array}\right\}
\quad\hbox{and}\quad
\left\{\begin{array}{l}
g_{n+1} := \alpha tg_n^2+h_n^2 \\[6pt]
h_{n+1} := \beta  tg_n^2+h_n^2
\end{array}\right\}
\quad \hbox{for all $n\ge0$.}$$

\begin{Lem}\label{MonicCase3}
Here $g_1=h_1=1$, and for every $n\ge2$ we have
\begin{eqnarray*}
g_n &=& \alpha^{2^{n-1}-1} \cdot t^{2^{n-1}-1} + \hbox{lower terms in $t$}, \\
h_n &=& \beta\cdot\alpha^{2^{n-1}-2} \cdot t^{2^{n-1}-1} + \hbox{lower terms in $t$}.
\end{eqnarray*}
\end{Lem}

\begin{Proof}
The recursion formula implies that $g_1=h_1=1$ and $g_2=\alpha t+1$ and $h_2=\beta t+1$. In particular the second statement holds for $n=2$. For arbitrary $n\ge2$ it follows by direct induction.
\end{Proof}

In projective coordinates the equation $f^k(0)=-f^\ell(0)$ now reads ${(g_k:h_k)} = {(-g_\ell:h_\ell)}$, which means that $g_k h_\ell + g_\ell h_k=0$. 
Here $\ell\ge1$, and we need one more case distinction:

\medskip
{\bf Case 3a:} $\ell\ge2$.\ \ 
Then from Lemma \ref{MonicCase3} and the fact that $k>\ell\ge2$ we deduce that 
$$g_k h_\ell + g_\ell h_k \ =\ 
2\cdot\beta\cdot\alpha^{2^{k-1}+2^{\ell-1}-3} \cdot t^{2^{k-1}+2^{\ell-1}-2} + \hbox{lower terms in $t$}.$$
With fixed constants $\alpha$, $\beta\not=0$ this equation has only finitely many solutions over a field of characteristic $\not=2$, as desired. 

\medskip
{\bf Case 3b:} $\ell=1$ and $\alpha+\beta\not=0$.\ \ 
Then $k>\ell=1$, and so Lemma \ref{MonicCase3} implies that 
$$g_k h_\ell + g_\ell h_k \ =\ g_k + h_k \ =\ 
(\alpha+\beta)\cdot\alpha^{2^{k-1}-2} \cdot t^{2^{k-1}-1} + \hbox{lower terms in $t$}.$$
For fixed constants $\alpha$, $\beta$ with $\alpha$, $\alpha+\beta\not=0$ this equation has only finitely many solutions over any field, as desired.

\medskip
{\bf Case 3c:} $\ell=1$ and $\alpha+\beta=0$.\ \ 
In this case we have $f(\infty) = \frac{\alpha}{\beta} = -1$ and hence $f^2(\infty) = f(-1) = f(1) = f^2(0)$. Since $k>\ell=1$, we can deduce from this and the equation $f^k(0) = - f^\ell(0)$ that $f^k(\infty) = f^k(0) = - f(0) = - 1 = f(\infty)$. By (\ref{RamLift}) this implies that $f^{k-1}(\infty) = \infty$ with $k-1>0$. 
%Thus the mapping scheme $\Gamma$ has the form
%$$\fbox{\ $\xymatrix@R-25pt@C=15pt{
%\scriptstyle i_0 & \scriptstyle i_1 &&& \\ 
%\bullet\ar[r]&\circbullet\ar[dddr]&&& \\
%{\phantom{x}}&&&&\ \ \ \\
%\scriptstyle j_0 & \scriptstyle j_1 & \ \scriptstyle i_2 && \\ 
%\bullet\ar[r]&\circbullet\ar[r]&\bullet\ar[r]&\ldots
%\ar[r]&\bullet\ar@/^17pt/[llll] \\
%\strut &&&&\\
%}$}\qquad\\$$
Setting $k' := 3(k-1)$ and $\ell' := 2(k-1)$, we then have $k'>\ell'\ge2$ and $f^{k'}(\infty) = \infty = - f^{\ell'}(\infty)$.
Thus after interchanging the roles of $0$ and $\infty$, which amounts to interchanging $i_0$ and $j_0$ in~$\Gamma$, we can replace $(k,\ell)$ by $(k',\ell')$. Afterwards we are back in the case 3a, which has already been completed.

\medskip
This finishes the proof of Theorem \ref{GammaModuliQuasiFinite}.

%%%%%%%%%%%%%%%%%%%%%%%%%%%%%%%%%%%%%%%%%%%%%%%%%%%%%%%%%%%%%%%%%%%%%%%%%%%%%%%%%%%%%%%%%%%%%
% \newpage

%%%%%%%%%%%%%%%%%%%%%%%%%%%%%%%%%%%%%%%%%%%%%%%%%%%%%%%%%%%%%%%%%%%%%%%%%%%%%%%%%%%%%%%%%%%%%

\end{document}